\RequirePackage{ifpdf}
\ifpdf % We are running pdfTeX in pdf mode
\documentclass[pdftex]{sigma}
\else
\documentclass{sigma}
\fi

\usepackage{rotate,eucal}
\usepackage{amsxtra}

\newcommand{\qdp}{\hspace*{1.5mm}}%

\newcommand{\qdn}{\hspace*{-1.5mm}}

\newcommand{\sss}{\scriptscriptstyle}

\newcommand{\mb}[1]{\mathbb{#1}}
\newcommand{\mc}[1]{\mathcal{#1}}

\newcommand{\fat}[3]{\left(#1;#2\right)_{#3}}
\newcommand{\ffat}[3]{\left[#1;#2\right]_{#3}}

\newcommand{\ffnk}[4]{\left[\qdn\ba{#1}#3\\#4\ea{\qdn\Big|#2}\right]}
\newcommand{\zbnm}[2]{\genfrac{(}{)}{0mm}{2}{#1}{#2}}

\newcommand{\zfrac}[2]{\genfrac{}{}{}{2}{#1}{#2}}
\newcommand{\ba}{\begin{array}}
\newcommand{\ea}{\end{array}}

\newcommand{\alp}{\alpha}
\newcommand{\bet}{\beta}
\newcommand{\gam}{\gamma}

\newcommand{\del}{\delta}

\newcommand{\vpi}{\varpi}

\newcommand{\Del}{\Delta}

\newtheorem{thm}{Theorem}%[section]

\newtheorem{corl}[thm]{Corollary}
\newtheorem{prop}[thm]{Proposition}

\newcommand{\bwd}{\bigtriangledown}

\begin{document}

\allowdisplaybreaks

\renewcommand{\PaperNumber}{050}

\FirstPageHeading

\ShortArticleName{Partial Sums of Two Quartic $q$-Series}

\ArticleName{Partial Sums of Two Quartic $\boldsymbol{q}$-Series}

\Author{Wenchang CHU~$^\dag$~and~Chenying WANG~$^\ddag$}

\AuthorNameForHeading{W. Chu and C. Wang}

\Address{$^\dag$~Dipartimento di Matematica, Universit\`a degli Studi di Salento,\\
\hphantom{$^\dag$}~Lecce-Arnesano P.~O.~Box~193, Lecce~73100, Italy}
\EmailD{\href{mailto:chu.wenchang@unile.it}{chu.wenchang@unile.it}}

\Address{$^\ddag$~College of Mathematics and Physics,  Nanjing University of Information Science\\
\hphantom{$^\ddag$}~and Technology, Nanjing 210044, P.~R.~China}
\EmailD{\href{mailto:wang.chenying@163.com}{wang.chenying@163.com}}

\ArticleDates{Received January 20, 2009, in f\/inal form April 17,
2009; Published online April 22, 2009}

\Abstract{The partial sums of two quartic basic hypergeometric
series are investigated by means of the modif\/ied Abel lemma
on summation by parts. Several summation and transformation
formulae are consequently established.}


\Classification{33D15; 05A15}

\section{Introduction and motivation}

For two indeterminate $x$ and $q$, the shifted-factorial of $x$
with base $q$ is def\/ined by
\[\fat{x}{q}{0} = 1
\qquad\text{and}\qquad
\fat{x}{q}{n} = (1-x) (1-xq)\cdots \big(1-xq^{n-1}\big)
\qquad\text{for}\quad
n\in\mb{N}.\]
When $|q|<1$, we have two well-def\/ined inf\/inite products
\[\fat{x}{q}{\infty}
 = \prod_{k=0}^{\infty}\big(1-q^kx\big)
\qquad\text{and}\qquad
\fat{x}{q}{n}
 = \fat{x}{q}{\infty}/\big(xq^n;q\big)_{\infty}.\]
The product and fraction of shifted factorials
are abbreviated respectively to
\begin{gather*}
\ffat{\alp,\bet,\dots,\gam}{q}{n}\qdp
 = \fat{\alp}{q}{n}\fat{\bet}{q}{n}\cdots\fat{\gam}{q}{n},\\
\ffnk{cccc}{q}{\alp, \bet, \dots, \gam}
          {A, B, \dots,C}_n
 =
\frac{\fat{\alp}{q}{n}\fat{\bet}{q}{n}\cdots\fat{\gam}{q}{n}}
     {\fat{A}{q}{n}\fat{B}{q}{n}\cdots\fat{C}{q}{n}}.
\end{gather*}
%%%%%%%%%%%%%%%%%%%%%%%%%%%%%%%%%%%%%%%%%%%%%%%%%%%%%%%%%%%%%%%%%
Following Gasper--Rahman~\cite{rahman}, the basic hypergeometric
series is def\/ined by
\[{_{1+r}\phi_s}
\ffnk{cccc}{q;z}{a_0, a_1, \dots, a_r}
                {b_1,\dots,b_s}
=\sum_{n=0}^\infty
\Big\{(-1)^nq^{\zbnm{n}{2}}\Big\}^{s-r}
\ffnk{cccc}{q}{a_0,a_1,\dots,a_r}
     {q,b_1,\dots,b_s}_nz^n, \]
where the base $q$ will be restricted to $|q|<1$
for nonterminating $q$-series. For its connections
to special functions and orthogonal polynomials,
the reader can refer, for example, to the monograph
written by Andrews--Askey--Roy~\cite{askey}
and the paper by Koornwinder~\cite{tom}.

In the theory of basic hypergeometric series, there are several
important classes, for example, well-poised~\cite{bailey},
quadratic~\cite{gessel83,rahman93v}, cubic~\cite{rahman93}
and quartic~\cite{gasper89,gasper90} series.
To our knowledge, there are four quartic series
which can be displayed as follows:
\begin{gather*}
F_n(a,b,d)
 := \sum_{k=0}^{n-1}\big(1-q^{5k}a\big)
\ffnk{c}{q}{b,d}{q^3a/b^2d^2}_k
\frac{[qa/bd;q]_{3k}}
     {[bd,bd/q,qbd;q^{2}]_k}
\ffnk{c}{q^4}{b^2d^2\!/q^2}{q^4a/b,q^4a/d}_k
q^{k};\\
G_n(a,c,e)
 := \sum_{k=0}^{n-1}\!\big(1-q^{5k}a\big)\!
\ffnk{c}{q}{c^2e^2\!/q^2a^3}{qa/c,qa/e}_k
\frac{[qa^2\!/ce,q^2a^2\!/ce,q^3a^2\!/ce;q^{2}]_k}
     {(ce/a;q)_{3k}}\!\!\!
\ffnk{c}{q^4}{c,e}{q^6a^4\!/c^2e^2}_k \! q^{k};\!\\
{U}_n(a,b,d)
 := \sum_{k=0}^{n-1}\!\big(1\!-\!q^{5k}a\big)\!
\frac{(q^2a/bd;q)_{k}}{(bd;q^2)_{2k}}\!
\ffnk{c}{q^2\!}{b,d}{q^5a^2/b^2d^2}_k\!\!\!
\frac{(q^3a^2/bd;q^6)_k}{(-bd/q^3a)^k}\!\!
\ffnk{c}{q^3\!}{b^2d^2/q^3a}{q^3a/b,q^3a/d}_k\!
\!q^{\zbnm{k}{2}};\!\\
{V}_n(a,c,e)
 := \sum_{k=0}^{n-1}\!\big(1-q^{5k}a\big)\!
\frac{(a^2/ce;q^2)_{2k}}{(qce/a;q)_k}
\ffnk{c}{q^2\!}{qc^2e^2/a^2}{qa/c,qa/e}_k
\frac{(-a/ce)^k}{(q^5ce;q^6)_k}\!
\ffnk{c}{q^3\!}{qc,qe}{q^2a^3/c^2e^2}_k
q^{-\zbnm{k}{2}}.
\end{gather*}
By means of the series rearrangement, Gasper and
Rahman~\cite{gasper89,gasper90,rahman}
discovered several summation and transformation formulae
for the nonterminating special cases of $F_n(a,b,d)$
and $U_n(a,b,d)$. The terminating series identities
for the last four sums have been established
by Chu~\cite{chu95b} and Chu--Wang~\cite{abelxy},
respectively, through inversion techniques and Abel's
lemma on summation by parts. For the partial sums
$F_n(a,b,d)$ and $G_n(a,c,e)$, the present authors~\cite{abelcy}
recently derived several useful reciprocal relations
and transformation formulae in terms of well-poised series.

The purpose of this paper is to investigate the remaining
two partial sums $U_n(a,b,d)$ and $V_n(a,c,e)$. By utilizing
the modif\/ied Abel lemma on summation by parts, we shall
show six unusual transformation formulae with two between
them and other four expressing $U_n(a,b,d)$ and $V_n(a,c,e)$
as partial sums of quadratic and cubic series. Several new
and known terminating as well as nonterminating series identities
are consequently obtained as particular instances.

In order to make the paper self-contained,
we reproduce Abel's lemma on summation by parts
(cf.~\cite{chu07a,chu08c,abelcy,abelxy})
as follows. For an arbitrary complex sequence $\{\tau_k\}$,
def\/ine the backward and forward dif\/ference operators $\bwd$
and $\Del$, respectively, by
\[
\bwd\tau_k=\tau_k-\tau_{k-1}
\qquad\text{and}\qquad
\Del\tau_k=\tau_{k+1}-\tau_k.
\]
Then \textit{Abel's lemma} on summation by parts can be
modif\/ied as follows:%\tag{$\star$}%%%\label{abelm}
\begin{gather*}
\sum_{k=0}^{n-1}B_k\bwd A_k
 =
\big\{{A_{n-1}B_n}-{A_{-1}B_0}\big\}
-\sum_{k=0}^{n-1}A_k\Del B_k.
\end{gather*}
This can be considered as the discrete
counterpart for the integral formula
\[
\int_a^bf(x)g'(x)\text{d}x
=\big\{f(b)g(b)-f(a)g(a)\big\}
-\int_a^bf'(x)g(x)\text{d}x.
\]

In fact, it is almost trivial to check the following expression
\[
\sum_{k=0}^{n-1}B_k\bwd A_k
 = \sum_{k=0}^{n-1}B_k\big\{A_{k}-A_{k-1}\big\}
 = \sum_{k=0}^{n-1}A_kB_{k}
 - \sum_{k=0}^{n-1}A_{k-1}B_k.
 \]
Replacing $k$ by $k+1$ for the last sum, we can reformulate
the equation as follows
\begin{gather*}
\sum_{k=0}^{n-1}B_k\bwd A_k
 = A_{n-1}B_{n}-A_{-1}B_0
+\sum_{k=0}^{n-1}A_k\big\{B_{k}-B_{k+1}\big\}\\
\phantom{\sum_{k=0}^{n-1}B_k\bwd A_k}{}
= A_{n-1}B_{n}-A_{-1}B_0-\sum_{k=0}^{n-1}A_k\Del B_k,
\end{gather*}
which is exactly the equality stated in the modif\/ied Abel lemma.

Throughout the paper, if $W_n$ is used to denote the
partial sum of some $q$-series, then the corresponding
letter $W$ without subscript will stand for the limit
of $W_n$ (if it exists of course) when $n\to\infty$.

When applying the modif\/ied Abel lemma on summation by parts
to deal with hypergeometric series, the crucial step lies
in f\/inding shifted factorial fractions $\{A_k,B_k\}$ so that
their dif\/ferences are expressible as ratios of linear factors.
This has not been an easy task, even though it is indeed routine
to factorize $\{A_k,B_k\}$ once they are f\/igured out.
Specif\/ically for  $U_n(a,b,d)$ and $V_n(a,b,d)$, we shall
devise three well-poised dif\/ference pairs for each partial sum.
This is based on numerous attempts to detect $A_k$ and $B_k$
sequences such that not only their dif\/ferences turn out to be
factorizable, but also their combinations match exactly the
summands displayed in $U_n(a,b,d)$ and $V_n(a,b,d)$.

The contents of the paper will be organized as follows.
In the second section, $U_n(a,b,d)$ will be reformulated
through the modif\/ied Abel lemma on summation by parts,
which lead to three transformations of $U_n(a,b,d)$ into
partial sums of quadratic, cubic and quartic series.
Then the third section will be devoted to the transformation
formulae of $V_n(a,b,d)$ in terms of partial sums of quadratic,
cubic and quartic series. These transformations on $U_n(a,b,d)$
and $V_n(a,b,d)$ will recover several known identities appeared
in Chu--Wang~\cite{chu95b,abelxy} and Gasper--Rahman~\cite{rahman},
and yield a few additional new summation formulae.

\section[Transformation and summation formulae for ${U}_n(a,b,d)$]{Transformation and summation formulae for $\boldsymbol{{U}_n(a,b,d)}$}\label{section2}

In this section, we shall investigate the partial
sum of quartic $q$-series ${U}_n(a,b,d)$. Three
transformation formulae will be established from
${U}_n(a,b,d)$ to quadratic, cubic and another
quartic series, unlike those for ${F}_n(a,b,d)$
and ${G}_n(a,c,e)$ shown in~\cite{abelcy}, where
reciprocal relations and transformations into
well-poised partial sums have been derived.
As particular cases of ${U}_n(a,b,d)$, four
nonterminating series will be evaluated, including
two that appeared in Gasper--Rahman's book~\cite{rahman}.

In order for the reader to gain an immediate insight
into the well-poised structure, we reformulate
the series ${U}_n(a,b,d)$ in the following manner:
\begin{gather*}
{U}_n(a,b,d)
 = \sum_{k=0}^{n-1}\big(1-q^{5k}a\big)
 \frac{[b,d;q^2]_k}
      {[q^3a/b,q^3a/d;q^3]_k}
\frac{(b^2d^2/q^3a;q^3)_k}
     {(q^5a^2\!/b^2d^2;q^2)_k}
q^{k}\\
\phantom{{U}_n(a,b,d)=}{} \times
\frac{[q^2a/bd,q^2a/bd;q]_k}
     {[bd,q^2bd;q^{4}]_k}
\frac{(q^3a^2/bd;q^6)_{k}}
     {(bd/q^2a;q^{-1})_k}.
\end{gather*}

\subsection{Quartic series to quadratic series} \label{section2.1}

Let $\text{A}_k$ and $\text{B}_k$ be def\/ined by
\begin{gather*}
\text{A}_k
 = \frac{[q^3a/bd,q^5a/bd;q]_k(b^3d^3/q^7a^2;q^2)_k(q^9a^2/bd;q^6)_k}
       {[bd,q^2bd;q^4]_k(q^{12}a^3/b^3d^3;q^3)_k(bd/q^4a;q^{-1})_k},\\
\text{B}_k
 = \frac{[b,d;q^2]_k}
        {[q^3a/b,q^3a/d;q^3]_k}
\frac{[b^2d^2/q^6a,q^{12}a^3/b^3d^3;q^3]_k}
     {[q^9a^2/b^2d^2,b^3d^3/q^9a^2;q^2]_k}.
\end{gather*}
We can easily show the relations
\begin{gather*}
\vpi
 := \text{A}_{-1}\text{B}_{0}
=\frac{a (1-q^2/bd)(1-q^4/bd)(1-q^3a/bd)(1-q^9a^3/b^3d^3)}
     {(1-q^2a/bd)(1-q^4a/bd)(1-q^3a^2/bd)(1-q^9a^2/b^3d^3)},\\
\mc{R}
 := \frac{\text{A}_{n-1}\text{B}_n}{\text{A}_{-1}\text{B}_0}
=\frac{1-q^{9+3n}a^3/b^3d^3}{1-q^9a^3/b^3d^3}
\ffnk{c}{q^2}{b, d}{q^9a^2/b^2d^2}_n\\
\phantom{\mc{R}  :=}{}  \times
\frac{[q^2a/bd,q^4a/bd;q]_n}
     {(bd/q^4;q^2)_{2n}}
\ffnk{c}{q^3}{b^2d^2/q^6a}{q^3a/b,q^3a/d}_n
\frac{(q^3a^2/bd;q^6)_n}{(bd/q^3a;q^{-1})_n};
\end{gather*}
and calculate the f\/inite dif\/ferences
\begin{gather*}
\bwd\text{A}_k
 =
\frac{(1-q^{5k}a)(1-q^{6-3k}a/b^2d^2)(1-q^{2k+5}a^2/b^2d^2)(1-q^{2k+7}a^2/b^2d^2)}
     {(1-q^2a/bd)(1-q^4a/bd)(1-q^3a^2/bd)(1-q^9a^2/b^3d^3)}\\
\phantom{\bwd\text{A}_k=}{}  \times
\frac{[q^2a/bd,q^4a/bd;q]_k(b^3d^3/q^9a^2;q^2)_k(q^3a^2/bd;q^6)_k}
     {[bd,q^2bd;q^4]_k(q^{12}a^3/b^3d^3;q^3)_k(bd/q^4a;q^{-1})_k}
q^{2k},\\
\Del\text{B}_k
 = -
\frac{(1-q^{3+5k}a)(1-q^{3+k}a/bd)(1-q^9a^2/b^2d^3)(1-q^9a^2/b^3d^2)}
     {(1-q^3a/b)(1-q^3a/d)(1-q^9a^2/b^2d^2)(1-q^9a^2/b^3d^3)}\\
 \phantom{\Del\text{B}_k=}{} \times
\ffnk{c}{q^2}{b,\qquad\:d}
               {q^{11}a^2/b^2d^2,b^3d^3/q^7a^2}_k
\ffnk{c}{q^3}{b^2d^2/q^6a,q^{12}a^3/b^3d^3}
             {q^6a/b,\quad\:q^6a/d}_kq^{2k}.
\end{gather*}
According to the modif\/ied Abel lemma on summation by parts,
the f\/inite $U$-sum can be reformulated as follows:
\begin{gather*}
U_n(a,b,d)
\frac{(1-q^5a^2/b^2d^2)(1-q^7a^2/b^2d^2)(1-q^6a/b^2d^2)}
     {(1-q^2a/bd)(1-q^4a/bd)(1-q^3a^2/bd)(1-q^9a^2/b^3d^3)}\\
\phantom{U_n(a,b,d)}{} = \sum_{k=0}^{n-1}\text{B}_k\bwd\text{A}_k
=\vpi\big\{\mc{R}-1\big\}
-\sum_{k=0}^{n-1}\text{A}_k\Del\text{B}_k.
\end{gather*}
Observing that there holds the equality for the last partial sum
\[
-\sum_{k=0}^{n-1}\text{A}_k\Del\text{B}_k
=\frac{{U}_n(q^3a,b,d)\:(1-q^3a/bd)(1-q^9a^2/b^2d^3)(1-q^9a^2/b^3d^2)}
      {(1-q^3a/b)(1-q^3a/d)(1-q^9a^2/b^2d^2)(1-q^9a^2/b^3d^3)},
      \]
we derive after some simplif\/ication the recurrence relation
\begin{gather*}
{U}_n(a,b,d)
 =
{U}_n(q^3a,b,d)
\frac{(q^2a/bd;q)_3(1-q^3a^2\!/bd)(1-q^9a^2\!/b^2d^3)(1-q^9a^2\!/b^3d^2)}
     {(q^5a^2\!/b^2d^2;q^2)_3(1-q^3a/b)(1-q^3a/d)(1-q^6a/b^2d^2)}
 \\
 \phantom{{U}_n(a,b,d)=}{} -a
\big\{1-\mc{R}(a,b,d)\big\}
\frac{(1-q^2/bd)(1-q^4/bd)(1-q^3a/bd)(1-q^9a^3/b^3d^3)}
     {(1-q^5a^2/b^2d^2)(1-q^7a^2/b^2d^2)(1-q^6a/b^2d^2)}.
\end{gather*}
Iterating this equation $m$-times, we get the following expression
\begin{gather*}
{U}_n(a,b,d)
 = {U}_n(q^{3m}a,b,d)
\frac{(q^2a/bd;q)_{3m}[q^3a^2\!/bd,q^9a^2\!/b^2d^3,q^9a^2\!/b^3d^2;q^6]_m}
     {(q^5a^2\!/b^2d^2;q^2)_{3m}[q^3a/b,q^3a/d,q^6a/b^2d^2;q^3]_m}\\
 \phantom{{U}_n(a,b,d)=}{} -
\frac{a\:(1-q^2/bd)(1-q^4/bd)(1-q^3a/bd)}
     {(1-q^5a^2/b^2d^2)(1-q^7a^2/b^2d^2)(1-q^6a/b^2d^2)}\\
 \phantom{{U}_n(a,b,d)=}{} \times \sum_{k=0}^{m-1}
(1-q^{9+9k}a^3\!/b^3d^3)
\ffnk{c}{q^3}{q^2a/bd,q^4a/bd,q^6a/bd}
     {q^3a/b,q^3a/d,q^9a/b^2d^2}_kq^{3k}\\
\phantom{{U}_n(a,b,d)=}{} \times \big\{1 - \mc{R}(q^{3k}a,b,d)\big\}
\ffnk{c}{q^6}{q^3a^2/bd,\:q^9a^2/b^2d^3,\:q^9a^2/b^3d^2}
     {q^9a^2/b^2d^2,q^{11}a^2/b^2d^2,q^{13}a^2/b^2d^2}_k.
\end{gather*}
Writing explicitly the $\mc{R}$-function
by separating $k$ and $n$ factorials
\begin{gather*}
\mc{R}(q^{3k}a,b,d)
 = \frac{1-q^{9+3n+9k}a^3/b^3d^3}{1-q^{9+9k}a^3/b^3d^3}
\ffnk{c}{q^2}{b,\:\:d}{q^{9+6k}a^2/b^2d^2}_n\\
\phantom{\mc{R}(q^{3k}a,b,d)=}{}  \times
\frac{[q^{2+3k}a/bd,q^{4+3k}a/bd;q]_n}
     {(bd/q^4;q^2)_{2n}}
\ffnk{c}{q^3}{q^{-6-3k}b^2d^2/a}{q^{3+3k}a/b,q^{3+3k}a/d}_n
\frac{(q^{3+6k}a^2/bd;q^6)_n}{(q^{-3-3k}bd/a;q^{-1})_n}\\
\phantom{\mc{R}(q^{3k}a,b,d)}{} = \frac{[q^2a/bd,q^4a/bd;q]_n}
       {(bd/q^4;q^2)_{2n}}
\ffnk{c}{q^2}{b,\:d}{q^9a^2/b^2d^2}_n
\ffnk{c}{q^3}{b^2d^2/q^6a}{q^3a/b,q^3a/d}_n
\frac{(q^3a^2/bd;q^6)_n}{(bd/q^3a;q^{-1})_n}\\
\phantom{\mc{R}(q^{3k}a,b,d)=}{}\times
\frac{1-q^{9+3n+9k}a^3\!/b^3d^3}{1-q^{9+9k}a^3\!/b^3d^3}
\ffnk{c}{q^3}{q^3a/b,q^3a/d,q^9a/b^2d^2}
     {q^2a/bd,q^4a/bd,q^6a/bd}_k
\frac{(q^{3+6n}a^2\!/bd;q^6)_k}{(q^3a^2\!/bd;q^6)_k}\\
\phantom{\mc{R}(q^{3k}a,b,d)=}{}\times
\frac{(q^9a^2/b^2d^2;q^2)_{3k}}{(q^{9+2n}a^2/b^2d^2;q^2)_{3k}}
\ffnk{c}{q^3}{q^{2+n}a/bd,q^{4+n}a/bd,q^{6+n}a/bd}
     {q^{3+3n}a/b,q^{3+3n}a/d,q^{9-3n}a/b^2d^2}_k,
\end{gather*}
and then def\/ining the partial sum of quadratic series in base $q^3$ by
\begin{gather*}
{U}_m^{\diamond}(a,b,d)
 = \sum_{k=0}^{m-1} (1-q^{9+9k}a^3\!/b^3d^3)
\ffnk{c}{q^3}{q^2a/bd,q^4a/bd,q^6a/bd}
     {q^3a/b,q^3a/d,q^9a/b^2d^2}_k\\
\phantom{{U}_m^{\diamond}(a,b,d)=}{} \times
\ffnk{c}{q^6}{q^3a^2/bd,q^9a^2/b^2d^3,q^9a^2/b^3d^2}
     {q^9a^2/b^2d^2,q^{11}a^2/b^2d^2,q^{13}a^2/b^2d^2}_k
q^{3k},
\end{gather*}
we derive the following transformation formula.

\begin{thm}[Transformation between quartic and quadratic series]\label{4u2}
\begin{gather*}
{U}_n(a,b,d)
 - {U}_n(q^{3m}a,b,d)
\frac{(q^2a/bd;q)_{3m}[q^3a^2\!/bd,q^9a^2\!/b^2d^3,q^9a^2\!/b^3d^2;q^6]_m}
     {(q^5a^2\!/b^2d^2;q^2)_{3m}[q^3a/b,q^3a/d,q^6a/b^2d^2;q^3]_m}\\
\phantom{{U}_n(a,b,d)}{}  =
\frac{(1-q^3a/bd)(1-bd/q^2)(1-bd/q^4)}
     {(1-q^5a^2\!/b^2d^2)(1-q^7a^2\!/b^2d^2)(1-b^2d^2/q^6a)}\\
\phantom{{U}_n(a,b,d)=}{} \times
\Bigg\{{U}_m^{\diamond}(a,b,d)-{U}_m^{\diamond}(q^{5n}a,q^{2n}b,q^{2n}d)
\frac{(q^3a^2/bd;q^6)_n}{(bd/q^3a;q^{-1})_n}\\
\phantom{{U}_n(a,b,d)=}{} \times
\ffnk{c}{q^2}{b,\:d}{q^9a^2/b^2d^2}_n
\frac{[q^2a/bd,q^4a/bd;q]_n}
     {(bd/q^4;q^2)_{2n}}
\ffnk{c}{q^3}{b^2d^2/q^6a}{q^3a/b,q^3a/d}_n\Bigg\}.
\end{gather*}
\end{thm}

%%%%%%%%%%%%%%%%%%%%%%%%%%%%%%%%%%%%%%%%%%%%%%%%%%%%%%%%%%%%%%%%%
By means of the Weierstrass $M$-test on uniformly convergent series
(cf.\ Stromberg~\cite[p.~141]{karl}), we can compute the following limit
\[
\lim_{m,n\to\infty}U_n(q^{3m}a,b,d)
= \sum_{k=0}^{\infty}(bd)^k
\frac{[b,d;q^2]_k}{(bd;q^2)_{2k}}
q^{4\zbnm{k}{2}}.\]
Letting $m,n\to\infty$ in Theorem~\ref{4u2}, we obtain
the transformation formula.
\begin{prop}[Nonterminating series transformation]
\begin{gather*}
{U}({a,b,d}) =
\frac{(1-q^3a/bd)(1-bd/q^2)(1-bd/q^4)}
     {(1-q^5a^2\!/b^2d^2)(1-q^7a^2\!/b^2d^2)(1-b^2d^2/q^6a)}
{U}^{\diamond}(a,b,d)\\
\phantom{{U}({a,b,d}) =}{}  +
\frac{(q^2a/bd;q)_{\infty}[q^3a^2\!/bd,q^9a^2\!/b^2d^3,q^9a^2\!/b^3d^2;q^6]_{\infty}}
     {(q^5a^2\!/b^2d^2;q^2)_{\infty}[q^3a/b,q^3a/d,q^6a/b^2d^2;q^3]_{\infty}}
\sum_{k=0}^{\infty}(bd)^k
\frac{[b,d;q^2]_k}{(bd;q^2)_{2k}}
q^{4\zbnm{k}{2}}.
\end{gather*}
\end{prop}

When $bd=q^{2+2\del}$ with $\del=0,1$, this proposition results in
\begin{gather*}
{U}({a,b,q^{2+2\del}/b})
 = \ffnk{c}{q^3}{a,\:\:q^{1-2\del}a}
               {q^3a/b,q^{1-2\del}ab}_{\infty}
\ffnk{c}{q^6}{q^{3-6\del}a^2b,q^{5-4\del}a^2/b}
             {q^{3-6\del}a^2, q^{5-4\del}a^2}_{\infty}\\
\phantom{{U}({a,b,q^{2+2\del}/b})=}{}  \times
\sum_{k=0}^{\infty}
\frac{[b,q^{2+2\del}/b;q^2]_k}{(q^{2+2\del};q^2)_{2k}}
q^{4\zbnm{k}{2}+(2+2\del)k}.
\end{gather*}
Taking $a=q^{1+4\del}$ in this equation and noting
that the initial condition
\[{U}({q^{1+4\del},b,q^{2+2\del}/b}) = 1-q^{1+4\del},\]
we recover the following formula due to Andrews~\cite[equation~(4.6)]{andrews}
and Ismail--Stanton~\cite[Proposition~6]{ismail}
\[
\sum_{k=0}^{\infty}
\frac{[b,q^{2+2\del}/b;q^2]_k}{(q^{2+2\del};q^2)_{2k}}
q^{4\zbnm{k}{2}+(2+2\del)k}
 = \ffnk{c}{q^6}{q^{2+2\del}b,q^{4+4\del}/b}
             {q^{2+2\del}, q^{4+4\del}}_{\infty},
             \]
which leads, under the replacement $a\to q^{2\del}a$,
to the nonterminating series formula.

\begin{corl}[Gasper--Rahman~\protect{\cite[Exercise~3.29(ii),(iii)]{rahman}}]
\label{rahman-x}
\begin{gather*}
\sum_{k=0}^{\infty}\frac{1\!-\!q^{5k+2\del}a}{1-q^{2\del}a}
\ffnk{c}{q^2}{b,q^{2+2\del}\!/b}{qa^2}_k
\frac{(a;q)_{k}(q^{1+2\del}\!/a;q^3)_k(q^{1+2\del}a^2;q^6)_{k}}
     {(q^{2+2\del};q^2)_{2k}[qab,q^{3+2\del}a/b;q^3]_k}
q^{\zbnm{k+1}{2}}(-a)^k\\
\qquad{}  =\ffnk{c}{q^3}{qa,\:q^{3+2\del}a}
               {qab,q^{3+2\del}a/b}_{\infty}
\ffnk{c}{q^6}{q^{5}a^2\!/b,q^{3-2\del}a^2b,q^{2+2\del}b,q^{4+4\del}\!/b}
             {q^{5}a^2,\:q^{3-2\del}a^2,\:q^{2+2\del},\:q^{4+4\del}}_{\infty}.
\end{gather*}
\end{corl}

\subsection{Quartic series to cubic series}

Def\/ine two sequences by
\begin{gather*}
\text{A}_k =
\frac{[q^3a/bd,bd^2/q^4a;q]_{k}(q^2b;q^2)_k(q^9a^2/bd;q^6)_k}
     {[q^2bd,q^9a^2/bd^2;q^4]_{k}(q^3a/b;q^3)_k(bd/q^4a;q^{-1})_k},\\
\text{B}_k =
\frac{(q^4a/bd;q)_k(d/q^2;q^2)_k(b^2d^2/q^3a;q^3)_k(q^9a^2/bd^2;q^4)_k}
     {(bd;q^4)_k(q^6a/d;q^3)_k(q^7a^2/b^2d^2;q^2)_k(bd^2/q^5a;q)_k}.
\end{gather*}
It is not hard to check the relations
\begin{gather*}
\vpi := \text{A}_{-1}\text{B}_{0}
=\frac{(1-bd/q^2)(1-q^5a^2/bd^2)(1-a/b)(1-bd/q^3a)}
      {(1-q^2a/bd)(1-bd^2/q^5a)(1-b)(1-q^3a^2/bd)},\\
\mc{R} :=
\frac{\text{A}_{n-1}\text{B}_{n}}{\text{A}_{-1}\text{B}_{0}}
=\frac{1-q^{5+4n}a^2/bd^2}{1-q^{5}a^2/bd^2}
\ffnk{c}{q^2}{b,\:d/q^2}{q^7a^2/b^2d^2}_n\\
\phantom{\mc{R} :=}{}  \times
\frac{[q^2a/bd,q^4a/bd;q]_n}
     {(bd/q^2;q^2)_{2n}}
\ffnk{c}{q^3}{b^2d^2/q^3a}{a/b,q^6a/d}_n
\frac{(q^3a^2/bd;q^6)_n}{(bd/q^3a;q^{-1})_n};
\end{gather*}
and compute the f\/inite dif\/ferences
\begin{gather*}
\bwd\text{A}_k =
\frac{(1-q^{5k}a)(1-q^{2-2k}/d)(1-q^{5+2k}a^2/b^2d^2)(1-q^{3+3k}a/d)}
     {(1-q^2a/bd)(1-q^5a/bd^2)(1-b)(1-q^3a^2/bd)}\\
\phantom{\bwd\text{A}_k = }{}  \times
\frac{[q^2a/bd,bd^2/q^5a;q]_{k}(b;q^2)_k(q^3a^2/bd;q^6)_k}
     {[q^2bd,q^9a^2/bd^2;q^4]_{k}(q^3a/b;q^3)_k(bd/q^4a;q^{-1})_k}
q^{k},\\
\Del\text{B}_k
 = -
\frac{(1-q^{4+5k}a)(1-q^{3+k}a/bd)(1-q^{2+2k}b)(1-q^{9}a^2/b^2d^3)}
     {(1-q^5a/bd^2)(1-q^7a^2/b^2d^2)(1-q^6a/d)(1-bd)}\\
\phantom{\Del\text{B}_k=}{} \times
\ffnk{c}{q}{q^4a/bd}{bd^2/q^4a}_k
\ffnk{c}{q^2}{d/q^2}{q^9a^2/b^2d^2}_k
\ffnk{c}{q^3}{b^2d^2/q^3a}{q^9a/d}_k
\ffnk{c}{q^4}{q^9a^2/bd^2}{q^4bd}_k
q^{k}.
\end{gather*}
Then applying the modif\/ied Abel lemma on summation by parts,
the ${U}$-sum can alternatively be reformulated as follows:
\begin{gather*}
{U}_n(a,b,d)
\frac{(1-q^2/d)(1-q^3a/d)(1-q^5a^2/b^2d^2)}
     {(1-q^2a/bd)(1-q^5a/bd^2)(1-b)(1-q^3a^2/bd)}\\
\phantom{{U}_n(a,b,d)}{}  = \sum_{k=0}^{n-1}\text{B}_k\bwd\text{A}_k
=\vpi\big\{\mc{R}-1\big\}
-\sum_{k=0}^{n-1}\text{A}_k\Del\text{B}_k.
\end{gather*}
%%%%%%%%%%%%%%%%%%%%%%%%%%%%%%%%%%%%%%%%%%%%%%%%%%%%%%%%%%%%%%%%%
Observing that the last partial sum results in
\[-\sum_{k=0}^{n-1}\text{A}_k\Del\text{B}_k
 = \frac{{U}_n(q^4a,q^4b,d/q^2)\:(1-q^2b)(1-q^3a/bd)(1-q^9a^2/b^2d^3)}
      {(1-bd)(1-q^6a/d)(1-q^5a/bd^2)(1-q^7a^2/b^2d^2)},\]
we derive after some simplif\/ication the relation
\begin{gather*}
{U}_n(a,b,d)
 = {U}_n(q^4a,q^4b,d/q^2)
\frac{(b;q^2)_2(q^2a/bd;q)_2(1-q^3a^2/bd)(1-q^9a^2/b^2d^3)}
     {(q^5a^2/b^2d^2;q^2)_2(q^3a/d;q^3)_2(1-q^2/d)(1-bd)}\\
\phantom{{U}_n(a,b,d)=}{}  + \big\{1-\mc{R}(a,b,d)\big\}
\frac{(1-bd/q^2)(1-q^5a^2/bd^2)(1-a/b)(1-q^3a/bd)}
     {(1-d/q^2)(1-q^3a/d)(1-q^5a^2/b^2d^2)}.
\end{gather*}
Iterating it $m$-times, we get the following expression
\begin{gather*}
{U}_n({a,b,d})
 = {U}_n({q^{4m}a,q^{4m}b,q^{-2m}d})
\frac{(b;q^2)_{2m}(q^2a/bd;q)_{2m}[q^3a^2\!/bd,q^9a^2\!/b^2d^3;q^6]_m}
     {(q^5a^2/b^2d^2;q^2)_{2m}(q^3a/d;q^3)_{2m}[q^2/d,bd;q^2]_m}\\
 \phantom{{U}_n({a,b,d})=}{} + \frac{(1-bd/q^2)(1-q^3a/bd)(1-a/b)}
       {(1-d/q^2)(1-q^3a/d)(1-q^5a^2/b^2d^2)}
\sum_{k=0}^{m-1}
\frac{(1-q^{5+8k}a^2/bd^2)(b;q^2)_{2k}}
     {(q^{7}a^2/b^2d^2;q^2)_{2k}}\\
\phantom{{U}_n({a,b,d})=}{} \times
\big\{1\!-\!\mc{R}(q^{4k}a,q^{4k}b,q^{\!-\!2k}d)\big\}\!\!
\ffnk{c}{q^2\!}{q^2a\!/\!bd,q^5a\!/\!bd}{q^4/d,\:bd/q^2}_k\!\!
\ffnk{c}{q^6\!}{q^3a^2\!/\!bd,q^9a^2\!/\!b^2\!d^3}{q^6a/d,\:q^9a/d}_k\!\!
q^{2k}.
\end{gather*}
%%%%%%%%%%%%%%%%%%%%%%%%%%%%%%%%%%%%%%%%%%%%%%%%%%%%%%%%%%%%%%%%%
Rewriting the $\mc{R}$-function explicitly as
\begin{gather*}
\mc{R}(q^{4k}a,q^{4k}b,q^{-2k}d)
 = \frac{1-q^{5+4n+8k}a^2/bd^2}{1-q^{5+8k}a^2/bd^2}
\ffnk{c}{q^2}{q^{4k}b,q^{-2-2k}d}{q^{7+4k}a^2/b^2d^2}_n\\
\phantom{\mc{R}(q^{4k}a,q^{4k}b,q^{-2k}d)=}{}  \times
\frac{[q^{2+2k}a/bd,q^{4+2k}a/bd;q]_n}
     {(q^{2k-2}bd;q^2)_{2n}}\!
\ffnk{c}{q^3}{b^2d^2/q^3a}{a/b,q^{6+6k}a/d}_n\!
\frac{(q^{3+6k}a^2/bd;q^6)_n}
     {(q^{-3-2k}bd/a;q^{-1})_n}\!\\
\phantom{\mc{R}(q^{4k}a,q^{4k}b,q^{-2k}d)}{} =
\frac{[q^2a/bd,q^4a/bd;q]_n}
     {(bd/q^2;q^2)_{2n}}\!
\ffnk{c}{q^2}{b,\:d/q^2}{q^7a^2/b^2d^2}_n\!
\ffnk{c}{q^3}{b^2d^2/q^3a}{a/b,q^6a/d}_n\!
\frac{(q^3a^2/bd;q^6)_n}{(bd/q^3a;q^{-1})_n}\!\\
\phantom{\mc{R}(q^{4k}a,q^{4k}b,q^{-2k}d)=}{} \times
\frac{1-q^{5+4n+8k}a^2/bd^2}{1-q^{5+8k}a^2/bd^2}
\ffnk{c}{q^2}{q^{2+n}a/bd,q^{5+n}a/bd,q^4/d,q^{-2}bd}
     {q^2a/bd,q^5a/bd,q^{4-2n}/d,q^{4n-2}bd}_k\\
\phantom{\mc{R}(q^{4k}a,q^{4k}b,q^{-2k}d)=}{} \times
\frac{(q^6a/d;q^3)_{2k}}
     {(q^{6+3n}a/d;q^3)_{2k}}
\ffnk{c}{q^2}{q^{2n}b,q^7a^2/b^2d^2}{b,q^{7+2n}a^2/b^2d^2}_{2k}
\frac{(q^{3+6n}a^2/bd;q^6)_k}{(q^{3}a^2/bd;q^6)_k},
\end{gather*}
and def\/ining further the f\/inite cubic sum in base $q^2$ by
\begin{gather*}
{U}_m^{\sss\triangle}(a,b,d)
 = \sum_{k=0}^{m-1} (1-q^{5+8k}a^2/bd^2)
\ffnk{c}{q^6}{q^3a^2/bd,q^9a^2/b^2d^3}{q^6a/d,\:q^9a/d}_k\\
\phantom{{U}_m^{\sss\triangle}(a,b,d)=}{} \times
\frac{(b;q^2)_{2k}}
     {(q^{7}a^2/b^2d^2;q^2)_{2k}}
\ffnk{c}{q^2}{q^2a/bd,q^5a/bd}{q^4/d,bd/q^2}_k
q^{2k},
\end{gather*}
we f\/ind the following transformation formula.
%%%%%%%%%%%%%%%%%%%%%%%%%%%%%%%%%%%%%%%%%%%%%%%%%%%%%%%%%%%%%%%%%%%%%%%
\begin{thm}[Transformation between quartic and cubic series]\label{4u3}
\begin{gather*}
{U}_n({a,b,d}) - {U}_n({q^{4m}a,q^{4m}b,q^{-2m}d})
\frac{(b;q^2)_{2m}(q^2a/bd;q)_{2m}[q^3a^2\!/bd,q^9a^2\!/b^2d^3;q^6]_m}
     {(q^5a^2\!/b^2d^2;q^2)_{2m}(q^3a/d;q^3)_{2m}[q^2\!/d,bd;q^2]_m}\\
\phantom{{U}_n({a,b,d})}{}  = \frac{(1-bd/q^2)(1-q^3a/bd)(1-a/b)}
       {(1\!-\!d/q^2)(1\!-\!q^3a/d)(1\!-\!q^5a^2\!/b^2d^2)}
\bigg\{{U}_m^{\sss\triangle}(a,b,d)
-{U}_m^{\sss\triangle}({q^{5n}a,q^{2n}b,q^{2n}d})\\
\phantom{{U}_n({a,b,d})=}{}\times
\frac{[q^2a/bd,q^4a/bd;q]_n}
     {(bd/q^2;q^2)_{2n}}
\ffnk{c}{q^2}{b,\:d/q^2}{q^7a^2/b^2d^2}_n
\ffnk{c}{q^3}{b^2d^2/q^3a}{a/b,q^6a/d}_n
\frac{(q^3a^2/bd;q^6)_n}
     {(bd/q^3a;q^{-1})_n}\bigg\}.
\end{gather*}
\end{thm}

%%%%%%%%%%%%%%%%%%%%%%%%%%%%%%%%%%%%%%%%%%%%%%%%%%%%%%%%%%%%%%%%%
By means of the Weierstrass $M$-test, we can compute the limit
\[\lim_{m,n\to\infty}U_n(q^{4m}a,q^{4m}b,q^{-2m}d)
 = \sum_{k=0}^{\infty}\Big(\frac{a}{b}\Big)^k
\frac{(b^2d^2/q^3a;q^3)_k}{(q^3a/b;q^3)_k}
q^{3\zbnm{k+1}{2}}.\]
The limiting case $m,n\to\infty$ of Theorem~\ref{4u3}
leads to the transformation formula.

\begin{prop}[Nonterminating series transformation]
\begin{gather*}
{U}({a,b,d}) =
\frac{(1-bd/q^2)(1-q^3a/bd)(1-a/b)}
       {(1-d/q^2)(1-q^3a/d)(1-q^5a^2/b^2d^2)}
{U}^{\sss\triangle}(a,b,d)\\
\phantom{{U}({a,b,d})}{}  +
\frac{(q^2a/bd;q)_{\infty}(b;q^2)_{\infty}[q^3a^2\!/bd,q^9a^2\!/b^2d^3;q^6]_{\infty}}
     {(q^3a/d;q^3)_{\infty}[q^2/d,bd,q^5a^2/b^2d^2;q^2]_{\infty}}
\sum_{k=0}^{\infty}\Big(\frac{q^3a}{b}\Big)^k
\frac{(b^2d^2\!/q^3a;q^3)_k}{(q^3a/b;q^3)_k}
q^{3\zbnm{k}{2}}.
\end{gather*}
\end{prop}

When $b=a$, this proposition reduces to the following relation
\begin{gather*}
{U}({a,a,d}) =
\frac{(q^2/d;q)_{\infty}(a;q^2)_{\infty}[q^3a/d,q^9/d^3;q^6]_{\infty}}
     {(q^3a/d;q^3)_{\infty}[q^2/d,ad,q^5/d^2;q^2]_{\infty}}
\sum_{k=0}^{\infty}
\frac{(ad^2\!/q^3;q^3)_k}{(q^3;q^3)_k}
q^{3\zbnm{k}{2}+3k}.
\end{gather*}
Taking $d=1$ in the last equation and noting ${U}({a,b,1})=1-a$, we get
\[\sum_{k=0}^{\infty}
\frac{(a/q^3;q^3)_k}{(q^3;q^3)_k}
q^{3\zbnm{k}{2}+3k}
 = \frac{(a;q^6)_{\infty}}{(q^3;q^6)_{\infty}},\]
which results also from  a limiting case of the
$q$-Bailey--Daum formula (cf.~\cite[II-9]{rahman}).

Combining the last two equations leads us to the
nonterminating series identity.

\begin{corl}[Gasper--Rahman~\protect{\cite[Exercise~3.29(i)]{rahman}}]\label{rahman-y}
\begin{gather*}
 \sum_{k=0}^{\infty}\frac{1-q^{5k}a}{1-a}
\ffnk{c}{q^2}{a,d}{q^5\!/d^2}_k
\ffnk{c}{q^3}{ad^2\!/q^3}{q^3,q^3a/d}_k
\frac{(q^2\!/d;q)_{k}(q^3a/d;q^6)_{k}}{(ad;q^2)_{2k}}
q^{\zbnm{k}{2}}\big(-{q^3\!/d}\big)^k\\
 \qquad{} =\ffnk{c}{q^2}{q^2a,q^3\!/d}{ad,q^5\!/d^2}_{\infty}
\ffnk{c}{q^6}{ad^2,q^9\!/d^3}{q^3,q^6a/d}_{\infty}.
\end{gather*}
\end{corl}

\subsection{Quartic series to quartic series}

Finally, for the two sequences given by
\begin{gather*}
\text{A}_k =
\frac{(q^3a/bd;q)_k(q^2b;q^2)_k(b^2d^2/a;q^3)_k(q^5a^2/b^2d;q^4)_k}
     {(q^2bd;q^4)_k(q^3a/b;q^3)_k(q^5a^2/b^2d^2;q^2)_k(b^2d/a;q)_k},\\
\text{B}_k =
\frac{[qa/bd,b^2d/a;q]_{k}(d/q^2;q^2)_k(q^3a^2/bd;q^6)_k}
     {[bd,qa^2/b^2d;q^4]_{k}(q^3a/d;q^3)_k(bd/q^2a;q^{-1})_k};
\end{gather*}
we have no dif\/f\/iculty to check the relations
\begin{gather*}
\vpi := \text{A}_{-1}\text{B}_{0}
=\frac{(1-b^2d/qa)(1-b^2d^2/q^3a^2)(1-b/a)(1-bd/q^2)}
      {(1-bd/q^2a)(1-b)(1-b^2d^2/q^3a)(1-b^2d/qa^2)},\\
\mc{R} :=
\frac{\text{A}_{n-1}\text{B}_{n}}{\text{A}_{-1}\text{B}_{0}}
=\frac{1-q^{n-1}b^2d/a}{1-q^{-1}b^2d/a}
\ffnk{c}{q^2}{b,\:d/q^2}{q^3a^2/b^2d^2}_n\\
\phantom{\mc{R} :=}{}  \times
\frac{[qa/bd,q^2a/bd;q]_n}
     {(bd/q^2;q^2)_{2n}}
\ffnk{c}{q^3}{b^2d^2/q^3a}{a/b,q^3a/d}_n
\frac{(q^3a^2/bd;q^6)_n}{(bd/q^2a;q^{-1})_n};
\end{gather*}
and compute the f\/inite dif\/ferences
\begin{gather*}
\bwd\text{A}_k =
\frac{(1-q^{5k}a)(1-q^{2-2k}/d)(1-q^{1+k}a/bd)(1-q^{3}a^2/b^3d^2)}
     {(1-b)(1-q^2a/bd)(1-qa^2/b^2d)(1-q^3a/b^2d^2)}\\
\phantom{\bwd\text{A}_k =}{}  \times
\ffnk{c}{q}{q^2a/bd}{b^2d/a}_k
\ffnk{c}{q^2}{b}{q^5a^2/b^2d^2}_k
\ffnk{c}{q^3}{b^2d^2/q^3a}{q^3a/b}_k
\ffnk{c}{q^4}{qa^2/b^2d}{q^2bd}_k
q^{3k},\\
\Del\text{B}_k
  = -
\frac{(1-q^{1+5k}a)(1-q^{3k}a/b)(1-q^{2+2k}b)(1-q^{3+2k}a^2/b^2d^2)}
     {(1-bd)(1-qa^2/b^2d)(1-q^3a/d)(1-q^2a/bd)}\\
\phantom{\Del\text{B}_k=}{}\times
\frac{[qa/bd,b^2d/a;q]_{k}(d/q^2;q^2)_k(q^3a^2/bd;q^6)_k}
     {[q^4bd,q^5a^2/b^2d;q^4]_{k}(q^6a/d;q^3)_k(bd/q^3a;q^{-1})_k}
q^{-k}.
\end{gather*}
Then by means of the modif\/ied Abel lemma on summation by parts,
the ${U}$-sum can be reformulated as follows:
\begin{gather*}
{U}_n(a,b,d)
\frac{(1-q^2/d)(1-qa/bd)(1-q^{3}a^2/b^3d^2)}
     {(1-b)(1-q^2a/bd)(1-qa^2/b^2d)(1-q^3a/b^2d^2)}\\
\phantom{{U}_n(a,b,d)}{}  = \sum_{k=0}^{n-1}\text{B}_k\bwd\text{A}_k
=\vpi\big\{\mc{R}-1\big\}
-\sum_{k=0}^{n-1}\text{A}_k\Del\text{B}_k.
\end{gather*}
Writing the last partial sum in terms of $U$-sum as
\[-\sum_{k=0}^{n-1}\text{A}_k\Del\text{B}_k
 =
\frac{{U}_n(qa,q^4b,d/q^2)\:(1-b/a)(1-q^{2}b)(1-b^2d^2/q^{3}a^2)}
     {(1-bd)(1-b^2d/qa^2)(1-q^3a/d)(1-bd/q^2a)},\]
we derive after some simplif\/ication the following relation
\begin{gather*}
{U}_n(a,b,d)
 = \big\{1-\mc{R}(a,b,d)\big\}
\frac{(1-a/b)(1-bd/q^2)(1-b^2d/qa)(1-b^2d^2/q^3a^2)}
     {(1-d/q^{2})(1-bd/qa)(1-b^3d^2/q^{3}a^2)}\\
\phantom{{U}_n(a,b,d)=}{}  - {U}_n(qa,q^4b,d/q^2)
\frac{({q^2a}/{bd})\:(b;q^2)_2(1-b^2d^2/q^3a)(1-b^2d^2/q^3a^2)(1-b/a)}
     {(1-bd)(1-q^2\!/\!d)(1-bd\!/\!qa)(1-q^3a\!/\!d)(1-b^3d^2\!/\!q^3a^2)}.
\end{gather*}
Iterating it $m$-times, we get the following expression
\begin{gather*}
{U}_n({a,b,d})
 = {U}_n({q^{m}a,q^{4m}b,q^{-2m}d})
\ffnk{c}{q^2}{b^2d^2\!/q^3a^2}{q^2/d,bd}_m\\
\phantom{{U}_n({a,b,d})=}{}  \times
\frac{(b;q^2)_{2m}}{[bd/qa,bd/q^2a;q]_m}
\ffnk{c}{q^3}{b/a,b^2d^2\!/q^3a}{q^3a/d}_m
\frac{(q^2a/bd;q^{-1})_m}{(b^3d^2/q^3a^2;q^6)_m}\\
\phantom{{U}_n({a,b,d})}{}  +
\frac{(1-a/b)(1-bd/q^2)(1-b^2d^2/q^3a^2)}
     {(1-d/q^{2})(1-bd/qa)(1-b^3d^2/q^{3}a^2)}
\sum_{k=0}^{m-1}(1-q^{5k-1}b^2d/a)
\ffnk{c}{q^2}{b^2d^2\!/qa^2}{q^4/d,bd/q^2}_k\\
\phantom{{U}_n({a,b,d})=}{} \times
\big\{1\!-\!\mc{R}(q^{k}a,q^{4k}b,q^{-2k}d)\big\}
\frac{(b;q^2)_{2k}}{(bd/a;q)_{k}}
\ffnk{c}{q^3}{q^3b/a,b^2d^2\!/q^3a}{q^3a/d}_k
\frac{\big(qa/bd\big)^k\:q^{-\zbnm{k}{2}}}
     {(q^3b^3d^2\!/a^2;q^6)_k}.
\end{gather*}
Separating $k$ and $n$ factorials in the $\mc{R}$-function
\begin{gather*}
\mc{R}(q^{k}a,q^{4k}b,q^{-2k}d)
 = \frac{1-q^{n+5k-1}b^2d/a}{1-q^{5k-1}b^2d/a}
\ffnk{c}{q^2}{q^{4k}b,q^{-2-2k}d}
             {q^{3-2k}a^2/b^2d^2}_n\\
\phantom{\mc{R}(q^{k}a,q^{4k}b,q^{-2k}d)=}{} \times
\frac{[q^{1-k}a/bd,q^{2-k}a/bd;q]_n}
     {(q^{2k-2}bd;q^2)_{2n}}\!
\ffnk{c}{q^3}{q^{3k-3}b^2d^2/a}
             {q^{-3k}a/b,q^{3+3k}a/d}_n\!
\frac{(q^3a^2/bd;q^6)_n}
     {(q^{k-2}bd/a;q^{-1})_n}\\
\phantom{\mc{R}(q^{k}a,q^{4k}b,q^{-2k}d)}{} =
\frac{[qa/bd,q^2a/bd;q]_n}
     {(bd/q^2;q^2)_{2n}}
\ffnk{c}{q^2}{b,\:d/q^2}{q^3a^2/b^2d^2}_n
\ffnk{c}{q^3}{b^2d^2/q^3a}{a/b,q^3a/d}_n
\frac{(q^3a^2/bd;q^6)_n}{(bd/q^2a;q^{-1})_n}\\
\phantom{\mc{R}(q^{k}a,q^{4k}b,q^{-2k}d)=}{} \times
\frac{1-q^{n+5k-1}b^2d/a}{1-q^{5k-1}b^2d/a}
\ffnk{c}{q^3}{q^{3n-3}b^2d^2/a,q^{3-3n}b/a,q^3a/d}
             {b^2d^2/q^3a,\:q^{3}b/a,\:q^{3+3n}a/d}_k
q^{nk}\\
\phantom{\mc{R}(q^{k}a,q^{4k}b,q^{-2k}d)=}{} \times
\frac{(bd/a;q)_{k}(q^{2n}b;q^2)_{2k}}
     {(q^{-n}bd/a;q)_{k}(b;q^2)_{2k}}
\ffnk{c}{q^2}{q^4/d,\:bd/q^2,\:q^{-1-2n}b^2d^2/a^2}
             {q^{4-2n}/d,q^{4n-2}bd,b^2d^2/qa^2}_k,
\end{gather*}
and then def\/ining the partial sum of quartic series
\begin{gather*}
{U}_m^{\star}(a,b,d)
 =\sum_{k=0}^{m-1}
(1-q^{5k-1}b^2d/a)
\ffnk{c}{q^2}{b^2d^2/qa^2}{q^4/d,bd/q^2}_k
\frac{(b;q^2)_{2k}}
     {(bd/a;q)_{k}}\\
\phantom{{U}_m^{\star}(a,b,d)=}{} \times
\ffnk{c}{q^3}{q^3b/a,b^2d^2/q^3a}{q^{3}a/d}_k
\frac{\big(-qa/bd\big)^k
 q^{-\zbnm{k}{2}}}{(q^3b^3d^2/a^2;q^6)_k},
\end{gather*}
we establish the following transformation formula.
%%%%%%%%%%%%%%%%%%%%%%%%%%%%%%%%%%%%%%%%%%%%%%%%%%%%%%%%%%%%%%%%%%%%%%%
\begin{thm}[Transformation between two quartic series]\label{4u4}
\begin{gather*}
{U}_n({a,b,d}) - {U}_n({q^{m}a,q^{4m}b,q^{\!-\!2m}d})
\ffnk{c}{q^2}{b^2d^2/q^3a^2}{bd,\:q^2/d}_m\\
\phantom{{U}_n({a,b,d})=}{}  \times
\frac{(b;q^2)_{2m}}
     {[bd/qa,bd/q^2a;q]_m}
\ffnk{c}{q^3}{b^2d^2/q^3a,b/a}{q^3a/d}_m
\frac{(q^2a/bd;q^{-1})_m}
     {(b^3d^2/q^3a^2;q^6)_m}\\
\phantom{{U}_n({a,b,d})}{} =
\frac{(1-a/b)(1\!-\!bd/q^2)(1\!-\!b^2d^2/q^3a^2)}
     {(1\!-\!d/q^{2})(1\!-\!bd/qa)(1\!-\!b^3d^2/q^{3}a^2)}
\bigg\{{U}_m^{\star}(a,b,d)
-{U}_m^{\star}({q^{5n}a,q^{2n}b,q^{2n}d})\\
\phantom{{U}_n({a,b,d})=}{} \times
\frac{[qa/bd,q^2a/bd;q]_n}
     {(bd/q^2;q^2)_{2n}}
\ffnk{c}{q^2}{b,\:d/q^2}{q^3a^2/b^2d^2}_n
\ffnk{c}{q^3}{b^2d^2/q^3a}{a/b,q^3a/d}_n
\frac{(q^3a^2/bd;q^6)_n}
     {(bd/q^2a;q^{-1})_n}\bigg\}.
\end{gather*}
\end{thm}

In particular for $m=n$, we have the reduced transformation.
%%%%%%%%%%%%%%%%%%%%%%%%%%%%%%%%%%%%%%%%%%%%%%%%%%%%%%%%%%%%%%%%%%%%%%%
\begin{prop}[Transformation between two quartic series: $\boldsymbol{q^3a^2=b^2d^2}$]\label{p:4u4}
\begin{gather*}
{U}_n({a,b,d})
 = {U}_n^{\star}({q^{5n}a,q^{2n}b,q^{2n}d})
\frac{(a;q^3)_{n}(bd;q^6)_n}
     {(q^3a/d;q^3)_{n}(q^3a/b;q^3)_{n-1}}\\
\phantom{{U}_n({a,b,d})=}{}  \times
\frac{(qa/bd;q)_n(q^3a/bd;q)_{n-1}}
     {(qa/bd;q^{-1})_n(bd;q^2)_{2n-1}}
\ffnk{c}{q^2}{q^2b,\:d}{q^2}_{n-1}.
\end{gather*}
\end{prop}

In order to examine the limiting case $n\to\infty$
of the last equation, we write
${U}_n^{\star}(q^{5n}a,q^{2n}b$, $q^{3/2+2n}a/b)$
explicitly as follows:
\begin{gather*}
{U}_n^{\star}({q^{5n}a,q^{2n}b,q^{3/2+2n}a/b})
 = \sum_{k=0}^{n-1}
(1-q^{n+5k+1/2}b)
\ffnk{c}{q^2}{q^{2-2n}}{q^{5/2-2n}b/a,q^{4n-1/2}a}_k\\
\phantom{{U}_n^{\star}({q^{5n}a,q^{2n}b,q^{3/2+2n}a/b})=}{}  \times
\ffnk{c}{q^3}{q^{3-3n}b/a,q^{3n}a}{q^{3/2+3n}b}_k
\frac{(q^{2n}b;q^2)_{2k}\big(-q^{n-1/2}\big)^k}
     {(q^{3/2-n};q)_{k}(q^6b;q^6)_k}
q^{-\zbnm{k}{2}}.
\end{gather*}
Inverting the summation index $k\to n-1-k$ and then applying
the relation
\begin{gather*}
\frac{(q^{2-2n};q^2)_{n-1-k}(q^{3-3n}b/a;q^3)_{n-1-k}}
     {(q^{5/2-2n}b/a;q^2)_{n-1-k}(q^{3/2-n};q)_{n-1-k}}
 \\
 \qquad{} = \frac{(q^{2};q^2)_{n-1}(q^{3}a/b;q^3)_{n-1}}
      {(q^{3/2}a/b;q^2)_{n-1}(q^{1/2};q)_{n-1}}
q^{1-n^2}
\frac{(q^{1/2};q)_{k}(q^{3/2}a/b;q^2)_{k}}
     {(q^{2};q^2)_{k}(q^{3}a/b;q^3)_{k}}
q^{k^2+2k},
\end{gather*}
we can reformulate the f\/inite sum ${U}_n^{\star}({q^{5n}a,q^{2n}b,q^{3/2+2n}a/b})$ as
\begin{gather*}
(-1)^{n-1}
\frac{(q^{2};q^2)_{n-1}(q^{3}a/b;q^3)_{n-1}}
     {(q^{3/2}a/b;q^2)_{n-1}(q^{1/2};q)_{n-1}}
q^{\frac{1-n^2}{2}}
\sum_{k=0}^{n-1}
(-1)^{k}
\frac{(q^{1\!/2};q)_{k}(q^{3/2}a\!/b;q^2)_{k}}
     {(q^{2};q^2)_{k}(q^{3}a/b;q^3)_{k}}
q^{\zfrac{k(k+2)}{2}}\\
\qquad{} \times
\frac{(1-q^{6n-5k-9/2}b)(q^{2n}b;q^2)_{2n-2-2k}(q^{3n}a;q^3)_{n-1-k}}
     {(q^{4n-1/2}a;q^2)_{n-1-k}(q^{3/2+3n}b;q^3)_{n-1-k}(q^6b;q^6)_{n-1-k}}.
     \end{gather*}
Substituting this expression into Proposition~\ref{p:4u4} and then
letting $n\to\infty$, we derive the following transformation formula
\begin{gather}
{U}({a,b,q^{3/2}a/b})
=
\frac{(b;q^2)_{\infty}(a;q^3)_{\infty}(q^{3/2}a;q^6)_{\infty}}
     {(q^{3/2}a;q^2)_{\infty}(q^{3/2}b;q^3)_{\infty}(b;q^6)_{\infty}}\nonumber\\
\phantom{{U}({a,b,q^{3/2}a/b})=}{} \times
\sum_{k=0}^{\infty}
\Big(\!-\!q^{3/2}\Big)^{k}
\frac{(q^{1/2};q)_{k}(q^{3/2}a/b;q^2)_{k}}
     {(q^{2};q^2)_{k}(q^{3}a/b;q^3)_{k}}
q^{\zbnm{k}{2}}.\label{star}
\end{gather}

From this transformation, we can derive two new interesting summation
formulae. First, taking $b=1$ in this equation and keeping in mind that
${U}({a,1,d})=1$, we obtain the following remarkable summation formula.

\begin{corl}[Nonterminating series identity]\label{Q-2F1(-1/3)}
\[\sum_{k=0}^{\infty}
\Big(-q^{3/2}\Big)^k
\frac{(q^{1/2};q)_k(q^{3/2}a;q^2)_k}
     {(q^2;q^2)_k(q^3a;q^3)_k}
q^{\zbnm{k}{2}}
=\frac{(q^{3/2}a;q^2)_{\infty}(q^{3/2};q^3)_{\infty}(q^6;q^6)_{\infty}}
      {(q^2;q^2)_{\infty}(q^3a;q^3)_{\infty}(q^{3/2}a;q^6)_{\infty}}.\]
\end{corl}
The special case $a=0$ of this corollary recovers an identity
of Rogers--Ramanujan type due to Stanton~\cite[p.~61]{stanton}:
\[\sum_{k=0}^{\infty}
\frac{(-q;q^2)_k}{(q^4;q^4)_k}q^{k(k+2)}
 = \frac{(-q;q^2)_{\infty}}{(q^2;q^2)_{\infty}}
[q^6,\:q,\:q^5;q^6]_{\infty}.\]
Combining \eqref{star} with Corollary~\ref{Q-2F1(-1/3)}
yields another formula for quartic series.

\begin{corl}[Nonterminating series identity]\label{QQ-2F1(-1/3)}
\begin{gather*}
 \sum_{k=0}^{\infty}(-1)^k
\frac{1-q^{5k}a}{1-a}
\ffnk{c}{q^2}{b,q^{3/2}a/b}{q^2}_k
\ffnk{c}{q^3}{a}{q^3a/b,q^{3/2}b}_k
\frac{(q^{1/2};q)_{k}(q^{3/2}a;q^6)_{k}}{(q^{3/2}a;q^2)_{2k}}
q^{\frac{k^2+2k}{2}}\\
\qquad{} =\ffnk{c}{q^2}{b,q^{3/2}a/b}{q^2,q^{3/2}a}_{\infty}
\ffnk{c}{q^3}{q^3a,q^{3/2}}{q^{3}a/b,q^{3/2}b}_{\infty}
\ffnk{c}{q^6}{q^6,q^{3/2}a}{b,q^{3/2}a/b}_{\infty}.
\end{gather*}
\end{corl}

\section[Transformation and summation formulae for ${V}_n(a,b,d)$]{Transformation and summation formulae for $\boldsymbol{{V}_n(a,b,d)}$}\label{section3}

The quartic series ${V}_n(a,c,e)$ may be considered
as dual one to the ${U}_n(a,b,d)$ in the last section
in the sense that the numerator factorials and
denominator factorials are inverted. This section
will be devoted analogously to investigation of
summation and transformation formulae for ${V}_n(a,c,e)$.
As the series ${U}_n(a,c,e)$, the following expression
for ${V}_n(a,c,e)$ makes its well-poised structure
more transparent
\begin{gather*}
{V}_n(a,c,e)
 = \sum_{k=0}^{n-1}(1-q^{5k}a)
 \frac{(qc^2e^2/a^2;q^2)_k}
      {(q^2a^3/c^2e^2;q^3)_k}
\frac{[qc,qe;q^3]_k}
     {[qa/c,qa/e;q^2]_k}
q^{k}\\
\phantom{{V}_n(a,c,e)=}{} \times
\frac{[a^2/ce,q^2a^2/ce;q^4]_{k}}
     {[qce/a,qce/a;q]_k}
\frac{(q^{-1}a/ce;q^{-1})_{k}}
     {(q^5ce;q^6)_k}.
\end{gather*}

\subsection{Quartic series to quadratic series}

Let $\text{A}_k$ and $\text{B}_k$ be def\/ined by
\begin{gather*}
\text{A}_k
 = \frac{[q^3c^2e^2/a^2,a^4/qc^3e^3;q^2]_k[q^4c,q^4e;q^3]_k}
        {[q^2a^3/c^2e^2,q^6c^3e^3/a^3;q^3]_k[qa/c,qa/e;q^2]_k},\\
\text{B}_k
 = \frac{(q^6c^3e^3/a^3;q^3)_k(a^2/ce;q^2)_{2k}(a/q^2ce;q^{-1})_k}
       {(a^4/q^3c^3e^3;q^2)_k[qce/a,q^3ce/a;q]_k(q^5ce;q^6)_k}.
\end{gather*}
We can easily show the following relations
\begin{gather*}
\vpi
 := \text{A}_{-1}\text{B}_{0}
=\frac{(1-a/qc)(1-a/qe)(1-a^3/qc^2e^2)(1-q^3c^3e^3/a^3)}
       {(1-qc)(1-qe)(1-qc^2e^2/a^2)(1-a^4/q^3c^3e^3)},\\
\mc{R}
 := \frac{\text{A}_{n-1}\mc{B}_{n}}{\text{A}_{-1}\mc{B}_{0}}
=\frac{1-q^{3+3n}c^3e^3/a^3}{1-q^{3}c^3e^3/a^3}
\ffnk{c}{q^2}{qc^2e^2/a^2}{a/qc,a/qe}_n\\
\phantom{\mc{R}  :=}{}  \times
\frac{(a^2/ce;q^2)_{2n}}
     {[qce/a,q^3ce/a;q]_n}
\ffnk{c}{q^3}{qc,qe}{a^3/qc^2e^2}_n
\frac{(a/q^2ce;q^{-1})_n}{(q^5ce;q^{6})_n};
\end{gather*}
and calculate the f\/inite dif\/ferences
\begin{gather*}
\bwd\text{A}_k
 = \frac{(1-q^{5k}a)(1-q^{2+k}ce/a)(1-q^2c^2e^3/a^3)(1-q^2c^3e^2/a^3)}
     {(1-qc)(1-qe)(1-qc^2e^2/a^2)(1-q^3c^3e^3/a^4)}\\
\phantom{\bwd\text{A}_k=}{}  \times
\ffnk{c}{q^2}{qc^2e^2/a^2,a^4/q^3c^3e^3}
               {qa/c,\quad\quad qa/e}_k
\ffnk{c}{q^3}{qc,\quad\quad qe}
             {q^2a^3/c^2e^2,q^6c^3e^3/a^3}_k q^{2k},\\
\Del\text{B}_k  = -
\frac{(1-q^{3+5k}a)(1-q^{1-3k}c^2e^2/a^3)
           (1-q^{3+2k}c^2e^2/a^2)(1-q^{5+2k}c^2e^2/a^2)}
     {(1-q^3c^3e^3/a^4)(1-qce/a)(1-q^3ce/a)(1-q^5ce)}\\
\phantom{\Del\text{B}_k  =}{} \times
\frac{(q^6c^3e^3/a^3;q^3)_k(a^2/ce;q^2)_{2k}(a/q^2ce;q^{-1})_k}
     {(a^4/qc^3e^3;q^2)_k[q^2ce/a,q^4ce/a;q]_k(q^{11}ce;q^6)_k}
q^{2k}.
\end{gather*}
Applying the modif\/ied Abel lemma on summation by parts,
we can manipulate the f\/inite $V$-sum as follows:
\begin{gather*}
{V}_n(a,c,e)
\frac{(1-q^{2}ce/a)(1-q^2c^2e^3/a^3)(1-q^2c^3e^2/a^3)}
     {(1-qc)(1-qe)(1-qc^2e^2/a^2)(1-q^3c^3e^3/a^4)}\\
\phantom{{V}_n(a,c,e)}{}  = \sum_{k=0}^{n-1}\text{B}_k\bwd\text{A}_k
=\vpi\big\{\mc{R}-1\big\}
-\sum_{k=0}^{n-1}\text{A}_k\Del\text{B}_k.
\end{gather*}
Noting that the last partial sum results in
\begin{gather*}
-
\sum_{k=0}^{n-1}\text{A}_k\Del\text{B}_k
\:=\:{V}_n(q^3a,q^3c,q^3e)
\frac{(1-qc^2e^2/a^3)(1-q^{3}c^2e^2/a^2)(1-q^{5}c^2e^2/a^2)}
     {(1-q^3c^3e^3/a^4)(1-qce/a)(1-q^3ce/a)(1-q^5ce)},
     \end{gather*}
we derive after some simplif\/ication the recurrence relation
\begin{gather*}
{V}_n(a,c,e)
 = {V}_n(q^3a,q^3c,q^3e)
\frac{(1-qc)(1-qe)(1-qc^2e^2/a^3)(qc^2e^2/a^2;q^2)_3}
     {(1-q^5ce)(1-q^2c^2e^3/a^3)(1-q^2c^3e^2/a^3)(qce/a;q)_3}\\
\phantom{{V}_n(a,c,e)=}{}  - \big\{1-\mc{R}(a,c,e)\big\}
\frac{a(1-qc/a)(1-qe/a)(1-qc^2e^2/a^3)(1-q^3c^3e^3/a^3)}
     {(1-q^{2}ce/a)(1-q^2c^3e^2/a^3)(1-q^2c^2e^3/a^3)}.
\end{gather*}
Iterating it $m$-times, we get the following expression
\begin{gather*}
{V}_n(a,c,e)
 = {V}_n(q^{3m}a,q^{3m}c,q^{3m}e)
\frac{[qc,qe,qc^2e^2/a^3;q^3]_{m}(qc^2e^2/a^2;q^2)_{3m}}
     {[q^5ce,q^2c^2e^3/a^3,q^2c^3e^2/a^3;q^6]_{m}(qce/a;q)_{3m}}\\
\phantom{{V}_n(a,c,e)=}{}  - \frac{a\:(1-qc/a)(1-qe/a)(1-qc^2e^2/a^3)}
     {(1-q^{2}ce/a)(1-q^2c^2e^3/a^3)(1-q^2c^3e^2/a^3)}\\
\phantom{{V}_n(a,c,e)=}{}  \times \sum_{k=0}^{m-1}
(1-q^{3+9k}c^3e^3/a^3)
\ffnk{c}{q^3}{qc,\:qe,\:q^4c^2e^2/a^3}
             {qce/a,q^3ce/a,q^5ce/a}_k
q^{3k}\\
\phantom{{V}_n(a,c,e)=}{} \times
\big\{1-\mc{R}(q^{3k}a,q^{3k}c,q^{3k}e)\big\}
\ffnk{c}{q^6}{qc^2e^2/a^2,q^3c^2e^2/a^2,q^5c^2e^2/a^2}
             {q^5ce,\:q^8c^2e^3/a^3,\:q^8c^3e^2/a^3}_k.
\end{gather*}
%%%%%%%%%%%%%%%%%%%%%%%%%%%%%%%%%%%%%%%%%%%%%%%%%%%%%%%%%%%%%%%%
Writing explicitly the $\mc{R}$-function by separating $k$ and
$n$ factorials
\begin{gather*}
\mc{R}(q^{3k}a,q^{3k}c,q^{3k}e)
 = \frac{1-q^{3+3n+9k}c^3e^3/a^3}{1-q^{3+9k}c^3e^3/a^3}
\ffnk{c}{q^2}{q^{1+6k}c^2e^2/a^2}{a/qc,\:a/qe}_n\\
\phantom{\mc{R}(q^{3k}a,q^{3k}c,q^{3k}e)=}{}  \times
\frac{(a^2/ce;q^2)_{2n}}
     {[q^{1+3k}ce/a,q^{3+3k}ce/a;q]_n}
\ffnk{c}{q^3}{q^{1+3k}c,q^{1+3k}e}
             {q^{-1-3k}a^3/c^2e^2}_n
\frac{(q^{-2-3k}a/ce;q^{-1})_n}
     {(q^{5+6k}ce;q^{6})_n}\\
\phantom{\mc{R}(q^{3k}a,q^{3k}c,q^{3k}e)}{} = \frac{(a^2/ce;q^2)_{2n}}
     {[qce/a,q^3ce/a;q]_n}
\ffnk{c}{q^2}{qc^2e^2/a^2}{a/qc,a/qe}_n
\ffnk{c}{q^3}{qc,qe}{a^3/qc^2e^2}_n
\frac{(a/q^2ce;q^{-1})_n}{(q^5ce;q^{6})_n}\\
\phantom{\mc{R}(q^{3k}a,q^{3k}c,q^{3k}e)=}{} \times
\frac{1-q^{3+3n+9k}c^3e^3\!/\!a^3}{1-q^{3+9k}c^3e^3/a^3}
\frac{(q^{1+2n}c^2e^2/a^2;q^2)_{3k}(q^{5}ce;q^6)_{k}}
     {(qc^2e^2/a^2;q^2)_{3k}(q^{5+6n}ce;q^6)_{k}}\\
\phantom{\mc{R}(q^{3k}a,q^{3k}c,q^{3k}e)=}{} \times
\ffnk{c}{q^3}{q^{1+3n}c,q^{1+3n}e,qce/a,q^3ce/a,q^5ce/a,q^{4-3n}c^2e^2/a^3}
             {qc,qe,q^{1+n}ce/a,q^{3+n}ce/a,q^{5+n}ce/a,q^{4}c^2e^2/a^3}_k,
\end{gather*}
and then def\/ining the f\/inite quadratic sum in base $q^3$ by
\begin{gather*}
{V}_m^{\diamond}(a,c,e)
 =\sum_{k=0}^{m-1}(1-q^{3+9k}c^3e^3/a^3)
\ffnk{c}{q^3}{qc,\:qe,\:q^4c^2e^2/a^3}
             {qce/a,q^3ce/a,q^5ce/a}_k\\
 \phantom{{V}_m^{\diamond}(a,c,e)=}{} \times
\ffnk{c}{q^6}{qc^2e^2/a^2,q^3c^2e^2/a^2,q^5c^2e^2/a^2}
             {q^5ce,\:q^8c^2e^3/a^3,\:q^8c^3e^2/a^3}_k
q^{3k},
\end{gather*}
we obtain the following transformation formula.

\begin{thm}[Transformation between quartic and quadratic series]\label{4v2}
\begin{gather*}
 {V}_n(a,c,e)
-{V}_n(q^{3m}a,q^{3m}c,q^{3m}e)
\frac{[qc,qe,qc^2e^2/a^3;q^3]_{m}(qc^2e^2/a^2;q^2)_{3m}}
     {[q^5ce,q^2c^2e^3/a^3,q^2c^3e^2/a^3;q^6]_{m}(qce/a;q)_{3m}}\\
\phantom{{V}_n(a,c,e)}{} =\frac{(1-a/qc)(1-a/qe)(1-qc^2e^2/a^3)}
     {(1\!-\!a/q^{2}ce)(1\!-\!q^2c^3e^2\!/a^3)(1\!-\!q^2c^2e^3\!/a^3)}
\Bigg\{{V}_m^{\diamond}(a,c,e)-{V}_m^{\diamond}(q^{5n}a,q^{3n}c,q^{3n}e)\\
\phantom{{V}_n(a,c,e)=}{}\times
\frac{(a^2/ce;q^2)_{2n}}
     {[qce/a,q^3ce/a;q]_n}
\ffnk{c}{q^2}{qc^2e^2/a^2}{a/qc,a/qe}_n
\ffnk{c}{q^3}{qc,qe}{a^3/qc^2e^2}_n
\frac{(a/q^2ce;q^{-1})_n}
     {(q^5ce;q^{6})_n}\Bigg\}.
\end{gather*}
\end{thm}

Letting $n\to 1+m$, $c\to a/q$ and $e\to q^{-1-3m}$ in this
theorem, we derive the summation formula, which does not seem
to have explicitly appeared previously.
%%%%%%%%%%%%%%%%%%%%%%%%%%%%%%%%%%%%%%%%%%%%%%%%%%%%%%%%%%%%%%%%
\begin{corl}[Terminating series identity]
\begin{gather*}
 \sum_{k=0}^{m}\frac{1-q^{5k}a}{1-a}
\ffnk{c}{q^2}{q^{-3-6m}}{q^{2},q^{2+3m}a}_k
\ffnk{c}{q^3}{a,q^{-3m}}{q^{6+6m}a}_k
\frac{(q^{2+3m}a;q^2)_{2k}}{(q^{-1-3m};q)_{k}}
\frac{(-1)^kq^{(2+3m)k-\zbnm{k}{2}}}
     {(q^{3-3m}a;q^6)_k}\\
\qquad{} =\ffnk{c}{q^3}{q^{6+3m}a,q^{-3m}/a}{q^2,\quad\:q^4}_m
\ffnk{c}{q^6}{q^5,\quad\:q^7}{q^{9+3m}a,q^{3-3m}/a}_m.
\end{gather*}
\end{corl}

\subsection{Quartic series to cubic series}
Def\/ine two sequences by
\begin{gather*}
\text{A}_k =
\frac{(a^2/qc^2e;q)_k(q^3c^2e^2/a^2;q^2)_k(q^4c;q^3)_k(q^4a^2/ce;q^4)_k}
     {(q^6c^2e/a;q^4)_k(q^2a^3/c^2e^2;q^3)_k(qa/c;q^2)_k(qce/a;q)_k},\\
\text{B}_k =
\frac{[q^2a^2/ce,q^6c^2e/a;q^4]_{k}(qe;q^3)_k(a/qce;q^{-1})_k}
     {[q^2ce/a,a^2/q^2c^2e;q]_{k}(q^3a/e;q^2)_k(q^5ce;q^6)_k}.
\end{gather*}
It is not hard to check the relations
\begin{gather*}
\vpi := \text{A}_{-1}\text{B}_{0}
=\frac{(1-ce/a)(1-a/qc)(1-a^3/qc^2e^2)(1-q^2c^2e/a)}
      {(1-a^2/q^2c^2e)(1-qc^2e^2/a^2)(1-qc)(1-a^2/ce)},\\
\mc{R} :=
\frac{\text{A}_{n-1}\text{B}_{n}}{\text{A}_{-1}\text{B}_{0}}
=\frac{1-q^{2+4n}c^2e/a}{1-q^2c^2e/a}
\ffnk{c}{q^2}{qc^2e^2/a^2}{a/qc,q^3a/e}_n\\
\phantom{\mc{R} :=}{} \times
\frac{(a^2/ce;q^2)_{2n}}
     {[ce/a,q^2ce/a;q]_n}
\ffnk{c}{q^3}{qc,\:qe}{a^3/qc^2e^2}_n
\frac{(a/qce;q^{-1})_n}{(q^5ce;q^6)_n};
\end{gather*}
and compute the f\/inite dif\/ferences
\begin{gather*}
\bwd\text{A}_k =
\frac{(1-q^{5k}a)(1-q^{1+k}ce/a)(1-q^{1+2k}a/e)(1-q^{2}c^3e^2/a^3)}
     {(1-q^2c^2e/a^2)(1-qc^2e^2/a^2)(1-qc)(1-a^2/ce)}\\
\phantom{\bwd\text{A}_k =}{} \times
\ffnk{c}{q}{a^2/q^2c^2e}{qce/a}_k
\ffnk{c}{q^2}{qc^2e^2/a^2}{qa/c}_k
\ffnk{c}{q^3}{qc}{q^2a^3/c^2e^2}_k
\ffnk{c}{q^4}{a^2/ce}{q^6c^2e/a}_k
q^{k},\\
\Del\text{B}_k
 = -
\frac{(1-q^{4+5k}a)(1-q^{1-2k}c/a)(1-q^{3+2k}c^2e^2/a^2)(1-q^{4+3k}c)}
     {(1-q^2ce/a)(1-q^2c^2e/a^2)(1-q^3a/e)(1-q^{5}ce)}\\
\phantom{\Del\text{B}_k=}{} \times
\frac{[q^2a^2/ce,q^6c^2e/a;q^4]_{k}(qe;q^3)_k(a/qce;q^{-1})_k}
     {[q^3ce/a,a^2/qc^2e;q]_{k}(q^5a/e;q^2)_k(q^{11}ce;q^6)_k}
q^{k}.
\end{gather*}
Then by means of the modif\/ied Abel lemma on summation by parts,
the ${V}$-sum can be reformulated as follows:
\begin{gather*}
{V}_n(a,c,e)
\frac{(1-qce/a)(1-qa/e)(1-q^{2}c^3e^2/a^3)}
     {(1-q^2c^2e/a^2)(1-qc^2e^2/a^2)(1-qc)(1-a^2/ce)}\\
\phantom{{V}_n(a,c,e)}{} = \sum_{k=0}^{n-1}\text{B}_k\bwd\text{A}_k
=\vpi\big\{\mc{R}-1\big\}
-\sum_{k=0}^{n-1}\text{A}_k\Del\text{B}_k.
\end{gather*}
By expressing the last partial sum in terms of $V$-sum as
\begin{gather*}
-
\sum_{k=0}^{n-1}\text{A}_k\Del\text{B}_k
\:=\:\frac{{V}_n(q^4a,q^6c,e)\:(1-qc/a)(1-q^{3}c^2e^2/a^2)(1-q^{4}c)}
        {(1-q^2ce/a)(1-q^2c^2e/a^2)(1-q^3a/e)(1-q^{5}ce)},
\end{gather*}
we derive after some simplif\/ication the following relation
\begin{gather*}
{V}_n(a,c,e)
 = {V}_n(q^4a,q^6c,e)
\frac{(qc;q^3)_2(qc^2e^2/a^2;q^2)_2(1-qc/a)(1-a^2/ce)}
     {(qce\!/a;q)_2(qa/e;q^2)_2(1-q^5ce)(1-q^2c^3e^2/a^3)}\\
\phantom{{V}_n(a,c,e)=}{}  + \big\{1-\mc{R}(a,c,e)\big\}
\frac{(1-ce/a)(1-qc/a)(1-a^3/qc^2e^2)(1-q^2c^2e/a)}
     {(1-a/qce)(1-qa/e)(1-q^2c^{3}e^2/a^3)}.
\end{gather*}
Iterating it $m$-times, we get the following expression
\begin{gather*}
{V}_n({a,c,e})
 = {V}_n({q^{4m}a,q^{6m}c,e})
\frac{(qc;q^{3})_{2m}(qc^2e^2/a^2;q^2)_{2m}[qc/a,a^2/ce;q^2]_m}
     {(qce/a;q)_{2m}(qa/e;q^2)_{2m}[q^5ce,q^2c^3e^2/a^3;q^6]_m}\\
\phantom{{V}_n({a,c,e})=}{}  +
\frac{(1-ce/a)(1-qc/a)(1-a^3\!/qc^2e^2)}
     {(1-a/qce)(1-qa/e)(1-q^2c^3e^2\!/a^3)}
\sum_{k=0}^{m-1}
\frac{(1\!-\!q^{2+8k}c^2e/a)(qc;q^{3})_{2k}}
     {[q^5ce,q^8c^3e^2\!/a^3;q^6]_k}\\
\phantom{{V}_n({a,c,e})=}{} \times \big\{1-\mc{R}(q^{4k}a,q^{6k}c,e)\big\}
\ffnk{c}{q^2}{q^3c/a,a^2/ce}{ce/a,q^3ce/a}_k
\frac{(qc^2e^2\!/a^2;q^2)_{2k}}{(q^3a/e;q^2)_{2k}}q^{2k}.
\end{gather*}
Rewriting the $\mc{R}$-function explicitly as
\begin{gather*}
\mc{R}(q^{4k}a,q^{6k}c,e)
 = \frac{1-q^{2+4n+8k}c^2e/a}{1-q^{2+8k}c^2e/a}
\ffnk{c}{q^2}{q^{1+4k}c^2e^2/a^2}
             {q^{-1-2k}a/c,q^{3+4k}a/e}_n\\
\phantom{\mc{R}(q^{4k}a,q^{6k}c,e)=}{} \times
\frac{(q^{2k}a^2/ce;q^2)_{2n}}
     {[q^{2k}ce/a,q^{2+2k}ce/a;q]_n}
\ffnk{c}{q^3}{q^{1+6k}c,qe}{a^3/qc^2e^2}_n
\frac{(q^{-1-2k}a/ce;q^{-1})_n}
     {(q^{5+6k}ce;q^6)_n}\\
\phantom{\mc{R}(q^{4k}a,q^{6k}c,e)}{} = \frac{(a^2/ce;q^2)_{2n}}
     {[ce/a,q^2ce/a;q]_n}
\ffnk{c}{q^2}{qc^2e^2/a^2}{a/qc,q^3a/e}_n
\ffnk{c}{q^3}{qc,\:qe}{a^3/qc^2e^2}_n
\frac{(a/qce;q^{-1})_n}{(q^5ce;q^6)_n}\\
\phantom{\mc{R}(q^{4k}a,q^{6k}c,e)=}{} \times
\frac{1-q^{2+4n+8k}c^2e/a}{1-q^{2+8k}c^2e/a}
\ffnk{c}{q^2}{ce/a,q^3ce/a,q^{4n}a^2/ce,q^{3-2n}c/a}
             {q^nce/a,q^{n+3}ce/a,a^2/ce,q^{3}c/a}_k\\
\phantom{\mc{R}(q^{4k}a,q^{6k}c,e)=}{}\times
\frac{(q^{1+3n}c;q^3)_{2k}(q^{5}ce;q^6)_{k}}
     {(qc;q^3)_{2k}(q^{5+6n}ce;q^6)_{k}}
\ffnk{c}{q^2}{q^{1+2n}c^2e^2/a^2,q^3a/e}
             {qc^2e^2/a^2,q^{3+2n}a/e}_{2k},
\end{gather*}
and def\/ining further the f\/inite cubic sum in base $q^2$ by
\begin{gather*} {V}_m^{\sss\triangle}(a,c,e)
=\sum_{k=0}^{m-1}
(1-q^{2+8k}c^2e/a)\!
\ffnk{c}{q^2}{q^3c/a,a^2/ce}{ce/a,q^3ce/a}_k
\frac{(qc;q^{3})_{2k}(qc^2e^2/a^2;q^2)_{2k}}
     {(q^3a/e;q^2)_{2k}[q^5ce,q^8c^3e^2/a^3;q^6]_k}
q^{2k},
\end{gather*}
we establish the following transformation formula.

\begin{thm}[Transformation between quartic and cubic series]\label{4v3}
\begin{gather*}
{V}_n({a,c,e})
 - {V}_n({q^{4m}a,q^{6m}c,e})
\frac{(qc;q^{3})_{2m}(qc^2e^2/a^2;q^2)_{2m}[qc/a,a^2/ce;q^2]_m}
     {(qce/a;q)_{2m}(qa/e;q^2)_{2m}[q^5ce,q^2c^3e^2/a^3;q^6]_m}\\
\phantom{{V}_n({a,c,e})}{}  =
\frac{(1-ce/a)(1-qc/a)(1-a^3\!/qc^2e^2)}
     {(1-a/qce)(1-qa/e)(1-q^2c^3e^2\!/a^3)}
\bigg\{{V}_m^{\sss\triangle}(a,c,e)
-{V}_m^{\sss\triangle}(q^{5n}a,q^{3n}c,q^{3n}e)\\
\phantom{{V}_n({a,c,e})=}{}  \times
\frac{(a^2/ce;q^2)_{2n}}{[ce/a,q^2ce/a;q]_n}
\ffnk{c}{q^2}{qc^2e^2/a^2}{a/qc,q^3a/e}_n
\ffnk{c}{q^3}{qc,qe}{a^3/qc^2e^2}_n
\frac{(a/qce;q^{-1})_n}{(q^5ce;q^6)_n}\bigg\}.
\end{gather*}
\end{thm}

When $n\to 1+\del+2m$, $e\to q^{-1/2}a^{3/2}/c$ and $c=q^{-1-3\del-6m}$
with $\del=0,1$, the last theorem recovers the following
summation formula.
%%%%%%%%%%%%%%%%%%%%%%%%%%%%%%%%%%%%%%%%%%%%%%%%%%%%%%%%%%%%%%%%
\begin{corl}[Chu~\protect{\cite[Equation~(4.8d)]{chu95b}}]
\begin{gather*}
 \sum_{k=0}^{n}\frac{1-q^{5k}a}{1-a}
\ffnk{c}{q^2}{a}{q^{2+3n}a,q^{\frac12-3n}/a^{\frac12}}_k
\frac{(q^{\frac12}a^{\frac12};q^2)_{2k}(-q^{\frac12}a^{-\frac12})^k}
     {(q^{\frac12}a^{\frac12};q)_{k}(q^{\frac92}a^{\frac32};q^6)_k}
\ffnk{c}{q^3}{q^{\frac32+3n}a^{\frac32},q^{-3n}}{q^{3}}_k
q^{-\zbnm{k}{2}}\\
\qquad{} =\begin{cases}
0,&n -\emph{odd};\\[2mm]
\ffnk{c}{q^2}{q^2a}{q^{3/2}a^{3/2}}_{3m}
\ffnk{c}{q^6}{q^3}{q^{9/2}a^{3/2}}_m,&n=2m.
\end{cases}
\end{gather*}
\end{corl}

%%%%%%%%%%%%%%%%%%%%%%%%%%%%%%%%%%%%%%%%%%%%%%%%%%%%%%%%%%%%%%%%%%%%%%%
Instead, letting $n\to 1+m$, $e\to q^{-1/2}a^{3/2}/c$ and $a=q^{-1-4m}$
in Theorem~\ref{4v3}, we recover another terminating series identity.
\begin{corl}[Chu and Wang~\protect{\cite[Corollary 40]{abelxy}}]
\begin{gather*}
 \sum_{k=0}^{m}\frac{1-q^{5k-1-4m}}{1-q^{-1-4m}}
\ffnk{c}{q^2}{q^{-1-4m}}{q^{-4m}/c,q^{2+2m}c}_k
\frac{(q^{-2m};q^2)_{2k}(-q^{1+2m})^k}
     {(q^{-2m};q)_{k}(q^{3-6m};q^6)_k}
\ffnk{c}{q^3}{qc,q^{-1-6m}}{q^{3}}_k
q^{-\zbnm{k}{2}}\\
\qquad{} =\ffnk{c}{q^6}{q^4c}{q^3}_m
\ffnk{c}{q^2}{q}{q^2c}_{2m}
\ffnk{c}{q^2}{q^2c}{q}_m.
\end{gather*}
\end{corl}

\subsection{Quartic series to quartic series}

Finally, for the two sequences given by
\begin{gather*}
\text{A}_k =
\frac{[q^4a^2/ce,q^2ce^2/a;q^4]_{k}(q^4c;q^3)_k(a/ce;q^{-1})_k}
     {[qce/a,q^3a^2/ce^2;q]_{k}(qa/c;q^2)_k(q^5ce;q^6)_k},\\
\text{B}_k =
\frac{(q^3a^2/ce^2;q)_k(c^2e^2/qa^2;q^2)_k(e/q^2;q^3)_k(q^2a^2/ce;q^4)_k}
     {(ce^2/q^2a;q^4)_k(q^2a^3/c^2e^2;q^3)_k(q^3a/e;q^2)_k(ce/qa;q)_k};
\end{gather*}
we have no dif\/f\/iculty to check the relations
\begin{gather*}
\vpi := \text{A}_{-1}\text{B}_{0}
=\frac{(1-ce/a)(1-q^2a^2/ce^2)(1-a/qc)(1-ce/q)}
      {(1-a^2/ce)(1-ce^2/q^2a)(1-qc)(1-qa/ce)},\\
\mc{R} :=
\frac{\text{A}_{n-1}\text{B}_{n}}{\text{A}_{-1}\text{B}_{0}}
=\frac{1-q^{2+n}a^2/ce^2}{1-q^{2}a^2/ce^2}
\ffnk{c}{q^2}{c^2e^2/qa^2}{a/qc,q^3a/e}_n\\
 \phantom{\mc{R} :=}{} \times
\frac{(a^2/ce;q^{2})_{2n}}
     {[ce/qa,ce/a;q]_{n}}
\ffnk{c}{q^3}{qc,\:e/q^2}{q^2a^3/c^2e^2}_n
\frac{(qa/ce;q^{-1})_n}{(ce/q;q^6)_n};
\end{gather*}
and compute the f\/inite dif\/ferences
\begin{gather*}
\bwd\text{A}_k =
\frac{(1-q^{5k}a)(1-q^{3k-2}e)(1-q^{1+2k}a/e)(1-q^{2k-1}c^2e^2/a^2)}
     {(1-a^2/ce)(1-ce^2/q^2a)(1-qc)(1-ce/qa)}\\
\phantom{\bwd\text{A}_k =}{}  \times
\frac{[a^2/ce,ce^2/q^2a;q^4]_{k}(qc;q^3)_k(qa/ce;q^{-1})_k}
     {[qce/a,q^3a^2/ce^2;q]_{k}(qa/c;q^2)_k(q^5ce;q^6)_k}
q^{-k},\\
\Del\text{B}_k = -
\frac{(1-q^{1+5k}a)(1-q^{k}ce/a)(1-q^{2k-1}a/c)(1-q^{4}a^3/c^2e^3)}
     {(1-qa/ce)(1-q^3a/e)(1-q^2a^3/c^2e^2)(1-ce^2/q^2a)}\\
\phantom{\Del\text{B}_k =}{} \times
\ffnk{c}{q}{q^3a^2/ce^2}{ce/a}_k
\ffnk{c}{q^2}{c^2e^2/qa^2}{q^5a/e}_k
\ffnk{c}{q^3}{e/q^2}{q^5a^3/c^2e^2}_k
\ffnk{c}{q^4}{q^2a^2/ce}{q^2ce^2/a}_k
q^{k}.
\end{gather*}
Then according to the modif\/ied Abel lemma on summation by parts,
the ${V}$-sum can be reformulated as follows:
\begin{gather*}
{V}_n(a,c,e)
\frac{(1-e/q^2)(1-qa/e)(1-c^2e^2/qa^2)}
     {(1-a^2/ce)(1-ce^2/q^2a)(1-qc)(1-ce/qa)}\\
 \phantom{{V}_n(a,c,e)}{} = \sum_{k=0}^{n-1}\text{B}_k\bwd\text{A}_k
=\vpi\big\{\mc{R}-1\big\}
-\sum_{k=0}^{n-1}\text{A}_k\Del\text{B}_k.
\end{gather*}
%%%%%%%%%%%%%%%%%%%%%%%%%%%%%%%%%%%%%%%%%%%%%%%%%%%%%%%%%%%%%%%%%
Observing that the last partial sum results in
\[-\sum_{k=0}^{n-1}\text{A}_k\Del\text{B}_k
 =
\frac{{V}_n(qa,q^3c,e/q^3)\:(1-ce/a)(1-a/qc)(1-q^{4}a^3/c^2e^3)}
     {(1-qa/ce)(1-q^3a/e)(1-q^2a^3/c^2e^2)(1-ce^2/q^2a)},\]
we derive after some simplif\/ication the following relation
\begin{gather*}
{V}_n(a,c,e)
 = a\frac{(1-a/ce)(1-q^2a^2/ce^2)(1-qc/a)(1-q/ce)}
      {(1-qa^2/c^2e^2)(1-q^2/e)(1-qa/e)}
\big\{\mc{R}(a,c,e)-1\big\}\\
\phantom{{V}_n(a,c,e)=}{}  - {V}_n(qa,q^3c,e/q^3)
\frac{(1-qc)(1-qc/a)(1-a/ce)(1-a^2/ce)(1-q^4a^3/c^2e^3)}
     {(ce/qa)(1-q^2a^3/c^2e^2)(1-qa^2/c^2e^2)(1-q^2/e)(qa/e;q^2)_2}.
\end{gather*}
Iterating it $m$-times, we get the following expression
\begin{gather*}
{V}_n({a,c,e})
 = {V}_n({q^{m}a,q^{3m}c,q^{-3m}e})
\ffnk{c}{q^2}{qc/a,a^2/ce}{qa^2/c^2e^2}_m\\
\phantom{{V}_n({a,c,e})=}{}  \times
\frac{[a/ce,qa/ce;q]_m}
     {(qa\!/\!e;q^2)_{2m}}
\ffnk{c}{q^3}{qc}{q^2a^3\!/\!c^2e^2,q^2\!/\!e}_n
\frac{(q^4a^3/c^2e^3;q^6)_m}{(ce/qa;q^{-1})_m}\\
\phantom{{V}_n({a,c,e})=}{} - a\frac{(1-a/ce)(1-qc/a)(1-q/ce)}
       {(1-qa^2/c^2e^2)(1-q^2/e)(1-qa/e)}\\
\phantom{{V}_n({a,c,e})=}{}\times\sum_{k=0}^{m-1}(1-q^{2+5k}a^2/ce^2)
\ffnk{c}{q^2}{q^3c\!/\!a,a^2\!/\!ce}{q^3a^2/c^2e^2}_k
\frac{(qa\!/\!ce;q)_{k}}
     {(q^3a/e;q^2)_{2k}}
q^{\zbnm{k}{2}}\\
\phantom{{V}_n({a,c,e})=}{} \times
\big\{1-\mc{R}(q^{k}a,q^{3k}c,q^{\!-\!3k}e)\big\}
\frac{(qc;q^3)_k(q^4a^3/c^2e^3;q^{6})_{k}}
     {[q^{2}a^3/c^2e^2,q^5/e;q^3]_k}
\Big(-\frac{q^2a}{ce}\Big)^k.
\end{gather*}
%%%%%%%%%%%%%%%%%%%%%%%%%%%%%%%%%%%%%%%%%%%%%%%%%%%%%%%%%%%%%%%%%
Writing explicitly the $\mc{R}$-function  as
\begin{gather*}
\mc{R}(q^{k}a,q^{3k}c,q^{-3k}e)
 = \frac{1-q^{2+n+5k}a^2/ce^2}{1-q^{2+5k}a^2/ce^2}
\ffnk{c}{q^2}{q^{-1-2k}c^2e^2/a^2}
             {q^{-1-2k}a/c,q^{3+4k}a/e}_n\\
\phantom{\mc{R}(q^{k}a,q^{3k}c,q^{-3k}e)=}{}  \times
\frac{(q^{2k}a^2/ce;q^{2})_{2n}}
     {[q^{-1-k}ce/a,q^{-k}ce/a;q]_{n}}
\ffnk{c}{q^3}{q^{1+3k}c,q^{-2-3k}e}
             {q^{2+3k}a^3/c^2e^2}_n
\frac{(q^{1+k}a/ce;q^{-1})_n}
     {(ce/q;q^6)_n}\\
\phantom{\mc{R}(q^{k}a,q^{3k}c,q^{-3k}e)}{}=
\frac{(a^2/ce;q^{2})_{2n}}
     {[ce/qa,ce/a;q]_{n}}
\ffnk{c}{q^2}{c^2e^2/qa^2}{a/qc,q^3a/e}_n
\ffnk{c}{q^3}{qc,\:e/q^2}{q^2a^3/c^2e^2}_n
\frac{(qa/ce;q^{-1})_n}{(ce/q;q^6)_n}\\
\phantom{\mc{R}(q^{k}a,q^{3k}c,q^{-3k}e)=}{}\times
\frac{1-q^{2+n+5k}a^2/ce^2}{1-q^{5k+2}a^2/ce^2}
\ffnk{c}{q^2}{q^{3}a^2/c^2e^2,q^{4n}a^2/ce,q^{3-2n}c/a}
     {q^{3-2n}a^2/c^2e^2,a^2/ce,q^{3}c/a}_k
q^{-nk}\\
\phantom{\mc{R}(q^{k}a,q^{3k}c,q^{-3k}e)=}{} \times
\frac{(q^{1-n}a/ce;q)_k(q^{3}a/e;q^2)_{2k}}
     {(qa/ce;q)_k(q^{3+2n}a/e;q^2)_{2k}}
\ffnk{c}{q^3}{q^{1+3n}c,q^{5}/e,q^2a^3/c^2e^2}
     {qc,q^{5-3n}/e,q^{2+3n}a^3/c^2e^2}_k,
\end{gather*}
and then def\/ining the f\/inite sum of quartic series
\begin{gather*}
{V}_m^{\star}(a,c,e)
 = \sum_{k=0}^{m-1}(1-q^{2+5k}a^2\!/ce^2)
\ffnk{c}{q^2}{q^3c/a,a^2/ce}{q^3a^2/c^2e^2}_k
\frac{(qa\!/\!ce;q)_{k}}
     {(q^3a/e;q^2)_{2k}}
q^{\zbnm{k}{2}}\\
\phantom{{V}_m^{\star}(a,c,e)=}{} \times
\Big(-\frac{q^2a}{ce}\Big)^k
\ffnk{c}{q^3}{qc}{q^{2}a^3/c^2e^2,q^5/e}_k
(q^4a^3/c^2e^3;q^{6})_{k},
\end{gather*}
we f\/ind the following transformation formula.

\begin{thm}[Transformation between two quartic series]\label{4v4}
\begin{gather*}
{V}_n({a,c,e}) - {V}_n({q^{m}a,q^{3m}c,q^{-3m}e})
\ffnk{c}{q^2}{qc/a,a^2/ce}{qa^2/c^2e^2}_m\\
\phantom{{V}_n({a,c,e}) -}{}  \times
\frac{[a/ce,qa/ce;q]_m}
     {(qa\!/\!e;q^2)_{2m}}
\ffnk{c}{q^3}{qc}{q^2a^3\!/\!c^2e^2,q^2\!/\!e}_n
\frac{(q^4a^3/c^2e^3;q^6)_m}{(ce/qa;q^{-1})_m}\\
\phantom{{V}_n({a,c,e})}{}=\frac{(1-a/qc)(1-a/ce)(1-ce/q)}
       {(1-qa^2\!/c^2e^2)(1-e\!/q^2)(1-qa\!/e)}
\bigg\{{V}_m^{\star}(a,c,e)
-{V}_m^{\star}({q^{5n}a,q^{3n}c,q^{3n}e})\\
\phantom{{V}_n({a,c,e})-}{}\times
\frac{(a^2/ce;q^{2})_{2n}}
     {[ce/qa,ce/a;q]_{n}}
\ffnk{c}{q^2}{c^2e^2/qa^2}{a/qc,q^3a/e}_n
\ffnk{c}{q^3}{qc,\:e/q^2}{q^2a^3/c^2e^2}_n
\frac{(qa/ce;q^{-1})_n}
     {(ce/q;q^6)_n}\bigg\}.
\end{gather*}
\end{thm}

Letting $n\to1+m$, $c=q^{-1-3m}$ and $e\to q^{2+3m}$, we obtain from
the last theorem the following terminating series identity.
%%%%%%%%%%%%%%%%%%%%%%%%%%%%%%%%%%%%%%%%%%%%%%%%%%%%%%%%%%%%%%%%
\begin{corl}[Terminating series identity]\label{nuova}
\begin{gather*}
 \sum_{k=0}^{m}\frac{1-q^{5k}a}{1-a}
\ffnk{c}{q^2}{q^{3}/a^2}{q^{2+3m}a,q^{-1-3m}a}_k
\ffnk{c}{q^3}{q^{3+3m},q^{-3m}}{a^3}_k
\frac{(a^2/q;q^2)_{2k}(-a)^k}
     {(q^{2}/a;q)_k(q^6;q^6)_k}
q^{-\zbnm{1+k}{2}}\\
\qquad{} =\frac{[a/q,qa;q]_{m}(q^{-3m}/a;q^2)_{m}(q^{-3m}a^3;q^6)_m}
         {(q^{-1-3m}a;q^2)_{2m}(a^3;q^3)_{m}(1/a;q^{-1})_{m}}.
\end{gather*}
\end{corl}

\section{Concluding remarks}

Recently, hypergeometric series has been found to have elliptic
analogue after the pioneering work of Frenkel--Turaev~\cite{turaev}.
Warnaar~\cite{warnaar} derived several terminating elliptic series
identities. Further summation formulae have been established
by Chu--Jia~\cite{chu08c} through Abel's lemma on summation
by parts. It is plausible that the same approach works also
for the elliptic analogue of the quartic series.
All what we have gotten are two terminating elliptic
analogues for Corollary~\ref{rahman-y} plus one
for Corollary~\ref{nuova}. Following the notations
of \cite{chu08c}, they are produced below
for reader's reference.

\begin{thm}[Terminating elliptic series identity]
\begin{gather*}
 \sum_{k=0}^{m}\frac{\theta(q^{5k}a;p)}{\theta(a;p)}
\ffnk{c}{q^2,p}{a,q^{2+m}}{q^{1-2m}}_k
\frac{[q^{-m};q,p]_k[q^{1-m}a;q^6,p]_k}
     {[q^{2+m}a;q^2,p]_{2k}}
q^{\zbnm{1+k}{2}-mk}\\
\qquad{} \times(-1)^k
\ffnk{c}{q^3,p}{q^{1+2m}a}{q^3,q^{1-m}a}_k
=
\frac{\chi(m=2n)[q^2a;q^2,p]_{n}[q^{3};q^{6},p]_{n}}
     {[q^{1+2n};q^2,p]_{n}[q^{4-2n}a;q^6,p]_n}.
\end{gather*}
\end{thm}
%%%%%%%%%%%%%%%%%%%%%%%%%%%%%%%%%%%%%%%%%%%%%%%%%%%%%%%%%%%%%%%%
\begin{thm}[Terminating elliptic series identity]
\begin{gather*}
 \sum_{k=0}^{m}\frac{\theta(q^{5k}a;p)}{\theta(a;p)}
\ffnk{c}{q^2,p}{a,q^{-2m}}{q^{5+4m}}_k
\frac{[q^{2+2m};q,p]_k[q^{3+2m}a;q^6,p]_k}
     {[q^{-2m}a;q^2,p]_{2k}}
q^{\zbnm{k}{2}+(3+2m)k}\\
\qquad{} \times (-1)^k
\ffnk{c}{q^3,p}{q^{-3-4m}a}{q^{3},q^{3+2m}a}_k
\:=\:
\frac{[q^2a;q^2,p]_{m}[q^{5};q^2,p]_{2m}[q^{-4m}a;q^{6},p]_{m}}
     {[q^{3};q^2,p]_{m}[q^{-2m}a;q^2,p]_{2m}[q^{9};q^6,p]_m}.
\end{gather*}
\end{thm}

\begin{thm}[Terminating elliptic series identity]
\begin{gather*}
 \sum_{k=0}^{m}\frac{\theta(q^{5k}a;p)}{\theta(a;p)}
\ffnk{c}{q^2,p}{q^{3}/a^2}{q^{2+3m}a,q^{-1-3m}a}_k
\frac{[a^2/q;q^2,p]_{2k}(-a)^k}
     {[q^{2}/a;q,p]_k[q^6;q^6,p]_k}
q^{-\zbnm{1+k}{2}}\\
\qquad{} \times
\ffnk{c}{q^3,p}{q^{3+3m},q^{-3m}}{a^3}_k
=\frac{[a/q,qa;q,p]_{m}[q^{-3m}/a;q^2,p]_{m}[q^{-3m}a^3;q^6,p]_m}
         {[q^{-1-3m}a;q^2,p]_{2m}[a^3;q^3,p]_{m}[1/a;q^{-1},p]_{m}}.
\end{gather*}
\end{thm}

\subsection*{Acknowledgments}

The work of the second author was partially supported by National Science Foundation of China (Youth grant 10801026).

\newpage

\pdfbookmark[1]{References}{ref}
\LastPageEnding

\end{document}